\documentclass[11pt]{article}

\usepackage{amsmath,amssymb,amsthm,bm,graphicx,color,epsfig}

\usepackage{algorithm}
\usepackage{algorithmic}

\usepackage{setspace}

\newtheorem{theorem}{Theorem}[section]

\newtheorem{algo}[theorem]{Algorithm}

\newcommand{\C}{\mathbb{C}}

\newcommand{\p}{\partial}
\newcommand{\Lapl}{\Delta}
\newcommand{\hf}{\frac{1}{2}}
\renewcommand{\P}{\mathcal{P}}
\newcommand{\G}{\mathcal{G}}

\newcommand{\wt}[1]{\widetilde{#1}}

\begin{document}

\title{Sweeping Preconditioner for the Helmholtz Equation:\\Moving
  Perfectly Matched Layers}
\author{Bj\"orn Engquist and Lexing Ying\\
  Department of Mathematics and ICES, University of Texas, Austin, TX
  78712}

\date{July 2010}
\maketitle

\begin{abstract}
  This paper introduces a new sweeping preconditioner for the
  iterative solution of the variable coefficient Helmholtz equation in
  two and three dimensions. The algorithms follow the general
  structure of constructing an approximate $LDL^t$ factorization by
  eliminating the unknowns layer by layer starting from an absorbing
  layer or boundary condition. The central idea of this paper is to
  approximate the Schur complement matrices of the factorization using
  moving perfectly matched layers (PMLs) introduced in the interior of
  the domain. Applying each Schur complement matrix is equivalent to
  solving a quasi-1D problem with a banded LU factorization in the 2D
  case and to solving a quasi-2D problem with a multifrontal method in
  the 3D case. The resulting preconditioner has linear application
  cost and the preconditioned iterative solver converges in a number
  of iterations that is essentially indefinite of the number of
  unknowns or the frequency. Numerical results are presented in both
  two and three dimensions to demonstrate the efficiency of this new
  preconditioner.
\end{abstract}



\section{Introduction}
\label{sec:intro}

This is the second of a series of papers on developing efficient
preconditioners for the numerical solutions of the Helmholtz equation
in two and three dimensions. To be specific, let the domain of
interest be the unit box $D = (0,1)^d$ with $d=2,3$. The
time-independent wave field $u(x)$ for $x\in D$ satisfies the
following Helmholtz equation,
\[
\Lapl u(x) + \frac{\omega^2}{c^2(x)} u(x) = f(x),
\]
where $\omega$ is the angular frequency, $c(x)$ is the velocity field
and, $f(x)$ is the external force. Commonly used boundary conditions
are the approximations of the Sommerfeld condition which guarantees
that the wave field generated by $f(x)$ propagates out of the domain
and other boundary condition for part of the boundary can also be
considered. By appropriately rescaling the system, it is convenient
to assume that the mean of $c(x)$ is around $1$. Then
$\frac{\omega}{2\pi}$ is the (average) wave number of this problem and
$\lambda = \frac{2\pi}{\omega}$ is the (typical) wavelength.

Equations of the Helmholtz type appear commonly in acoustics,
elasticity, electromagnetics, geophysics, and quantum
mechanics. Efficient and accurate numerical solution of the Helmholtz
equation is a very important problem in current numerical
mathematics. This is, however, a very difficult computational task due
to two main reasons. First, in a typical setting, the Helmholtz
equation is discretized with at least a constant number of points per
wavelength. Therefore, the number of samples $n$ in each dimension is
proportional to $\omega$, the total number of samples $N$ is $n^d =
O(\omega^d)$, and the approximating discrete system of the Helmholtz
equation is an $O(\omega^d) \times O(\omega^d)$ linear system, which
is extremely large in many practical high frequency
simulations. Second, since the discrete system is highly indefinite
and has a very oscillatory Green's function due to the wave nature of
the Helmholtz equation, most direct and iterative solvers developed
based on the multiscale paradigm are no longer efficient anymore. For
further remarks, see the discussion in \cite{EngquistYing:10a}.

\subsection{Approach and contribution}

In the previous paper \cite{EngquistYing:10a}, we introduced a sweeping
preconditioner that constructs an approximate $LDL^t$ factorization
layer by layer starting from an absorbing layer. An important
observation regarding the sweeping preconditioner is that the
intermediate Schur complement matrices of the $LDL^t$ factorization
corresponds to the restriction of the half-space Green's function of
the Helmholtz equation to a single layer. In \cite{EngquistYing:10a},
we represented the intermediate Schur complement matrices of the
factorization efficiently in the hierarchical matrix framework
\cite{Hackbusch:99}. In 2D, the efficiency of this preconditioner is
supported by analysis, has linear complexity, and results very small
number of iterations when combined with the GMRES solver. In 3D,
however, the theoretical justification is lacking and constructing the
preconditioner can be more costly.

In this paper, we propose a new sweeping preconditioner that works
well in both two and three dimensions. The central idea of this new
approach is to represent these Schur complement matrices in terms of
{\em moving perfectly matched layers} introduced in the interior of
the domain. Applying these Schur complement matrices then corresponds
to inverting a discrete Helmholtz system of a moving PML. Since each
moving PML is only of a few grids wide, fast direct algorithms can be
leveraged for this task. In 2D, this discrete system of the moving PML
layer is a quasi-1D problem and can be solved efficiently using a
banded LU factorization in an appropriate ordering. The construction
and application costs of the preconditioner are $O(n^2) = O(N)$ and
$O(n^2) = O(N)$, respectively. In 3D, the discrete Helmholtz system of
the moving PML is a quasi-2D problem and can be solved efficiently
using the multifrontal methods. The construction and application costs
of the preconditioner are $O(n^4) = O(N^{4/3})$ and $O(n^3 \log n) =
O(N\log N)$, respectively. Numerical results show that in both 2D and
3D this new sweeping preconditioner gives rise to iteration numbers
that is essentially independent of $N$ when combined with the GMRES
solver. After the construction of the preconditioner, we thus have a
linear solution method for the discrete Helmholtz system.

\subsection{Related work}

There has been a vast literature on developing efficient algorithms
for the Helmholtz equation. A partial list of significant progresses
includes
\cite{BaylissGoldsteinTurkel:83,BenamouDespres:97,BrandtLivshits:97,Despres:91,ElmanErnstOLeary:01,
  EngquistYing:07,ErlanggaOosterleeVuik:06,LairdGiles:02,
  OseiKuffuorSaad:09,Rokhlin:93}.  We refer to
the review article \cite{Erlangga:08} and our previous paper
\cite{EngquistYing:10a} for detailed discussion. The brief discussion
below is restricted to the ones that are closely related to the
approach proposed in this paper.

The most efficient direct methods for solving the discrete Helmholtz
systems are the multifrontal methods or their pivoted versions
\cite{DuffReid:83,George:73,Liu:92}. The multifrontal methods exploit
the locality of the discrete operator and construct an $LDL^t$
factorization based on a hierarchical partitioning of the domain. The
cost of a multifrontal method depends strongly on the number of
dimensions. For a 2D problem with $N=n^2$ unknowns, a multifrontal
method takes $O(N^{3/2})$ flops and $O(N \log N)$ storage space. The
prefactor is usually rather small, making the multifrontal methods
effectively the default choice for most 2D Helmholtz
problems. However, for a 3D problem with $N=n^3$ unknowns, a
multifrontal method requires $O(N^2)$ flops and $O(N^{4/3})$ storage
space, which can be very costly for large scale 3D problems.

The approach proposed here essentially reduces the dimensions of the
problem by working with $n$ subproblems with one dimension lower.  In
the 3D case, for each subproblem, it leverages the effectiveness of
the 2D multifrontal methods by solving a quasi-2D problem. The price
of this reduction is that we only end up with an approximate
inverse. However, this approximate inverse is reasonably accurate and
works very well as a preconditioner when combined with standard
iterative solvers in all our variable coefficient test cases.

\subsection{Contents}

The rest of this paper is organized as follows. Section
\ref{sec:2Dpre} presents the new sweeping preconditioner in the 2D
case and Section \ref{sec:2Dnum} reports the 2D numerical results. We
extend this approach to the 3D case in Section \ref{sec:3Dpre} and
report the 3D numerical results in Section \ref{sec:3Dnum}. Finally,
Section \ref{sec:conc} discusses some future directions of this work.

\section{Preconditioner in 2D}
\label{sec:2Dpre}

We will first discuss the sweeping factorization in general and then
introduce the moving PML.

\subsection{Discretization and sweeping factorization}

\begin{figure}[h!]
  \begin{center}
    \includegraphics{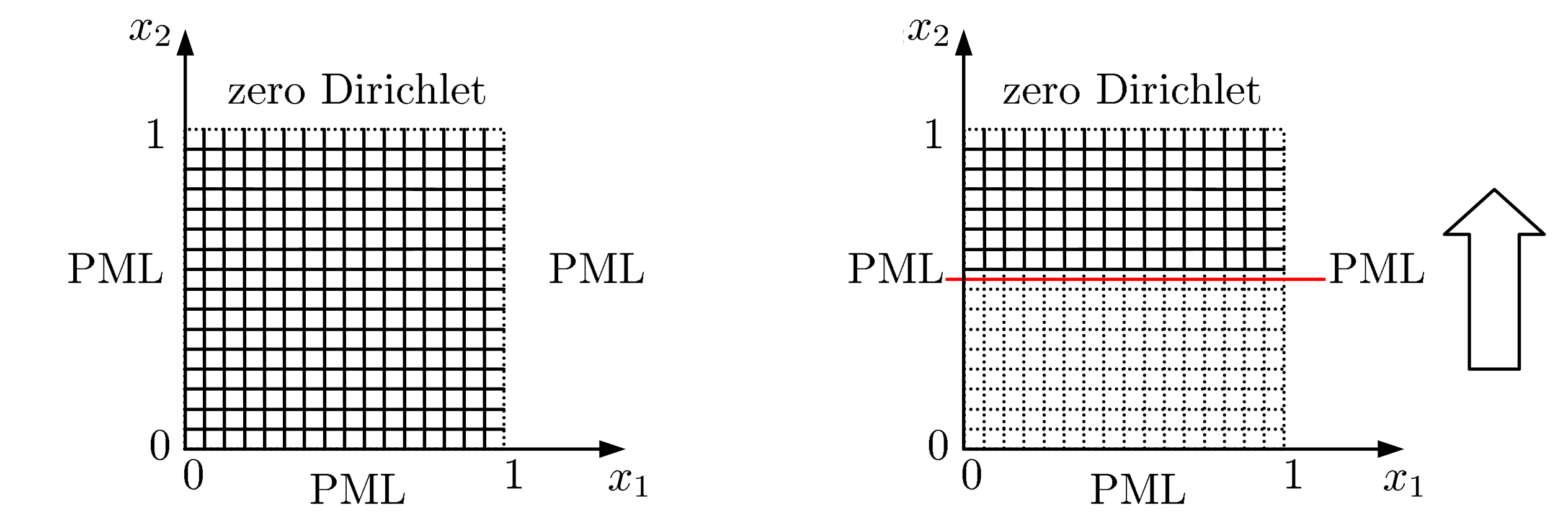}
  \end{center}
  \caption{Left: Discretization grid in 2D. Right: Sweeping order in
    2D with the moving PML. The dotted grid indicates the part that
    has already been eliminated.}
  \label{fig:2Dgrid}
\end{figure}

Recall that our computational domain in 2D is $D=(0,1)^2$. In order to
simplify the discussion, we assume that the Dirichlet zero boundary
condition is used on the side $x_2=1$ while approximations to the
Sommerfeld boundary condition is enforced on the other three
sides. One standard way of incorporating the Sommerfeld boundary
condition is to use the perfectly matched layer (PML)
\cite{Berenger:94,ChewWeedon:94,Johnson:10}. Introduce
\begin{equation}
\sigma_1(t) = 
\begin{cases}
  \frac{C}{\eta}\cdot \left( \frac{t-\eta}{\eta} \right)^2 & t \in [0,\eta]\\
  0 & t\in [\eta, 1-\eta] \\
  \frac{C}{\eta}\cdot \left( \frac{t-1+\eta}{\eta} \right)^2 & t \in [1-\eta,1]
\end{cases},
\quad
\sigma_2(t) = 
\begin{cases}
  \frac{C}{\eta}\cdot \left( \frac{t-\eta}{\eta} \right)^2 & t \in [0,\eta]\\
  0 & t\in [\eta, 1]
\end{cases},
\label{eq:sigma}
\end{equation}
and 
\[
s_1(x_1) = \left( 1+i\frac{\sigma_1(x_1)}{\omega} \right)^{-1},\quad
s_2(x_2) = \left( 1+i\frac{\sigma_2(x_2)}{\omega} \right)^{-1}.
\]
Here $\eta$ is typically around one wavelength and $C$ is an
appropriate positive constant independent of $\omega$. The PML method
replaces $\p_1$ with $s_1(x_1) \p_1$ and $\p_2$ with $s_2(x_2)\p_2$,
respectively. This effectively provides a damping layer of width
$\eta$ near the three sides with the Sommerfeld boundary
condition. The resulting equation becomes
\begin{eqnarray*}
\left( (s_1\p_1)(s_1\p_1) + (s_2\p_2)(s_2\p_2) + \frac{\omega^2}{c^2(x)} \right) u = f && x\in D=(0,1)^2,\\
u = 0 && x \in \p D.
\end{eqnarray*}
We assume that $f(x)$ is supported inside
$[\eta,1-\eta]\times[\eta,1]$ (away from the PML). Dividing the above
equation by $s_1(x_1) s_2(x_2)$ results
\[
\left( \p_1\left(\frac{s_1}{s_2} \p_1\right) + \p_2\left(\frac{s_2}{s_1} \p_2\right) +   \frac{\omega^2}{s_1s_2c^2(x)} \right) u = f.
\]
The main advantage of this equation is its symmetry. We discretize the
domain $[0,1]^2$ with a Cartesian grid with spacing $h = 1/(n+1)$. The
number of points $n$ in each dimension is proportional to the wave
number $\omega$ since a constant number of points is required for each
wavelength. The set of all interior points of this grid is denoted by 
\[
\P = \{ p_{i,j}=(ih,jh): 1\le i,j \le n\}
\]
(see Figure \ref{fig:2Dgrid} (left)) and the total number of grid
points is $N=n^2$.

We denote by $u_{i,j}$, $f_{i,j}$, and $c_{i,j}$ the values of $u(x)$,
$f(x)$, and $c(x)$ at point $p_{i,j} = (ih,jh)$. The $5$-point stencil
finite difference method writes down the equation at points in $\P$
using central difference. The resulting equation at $x_{i,j}=(ih,jh)$
is
\begin{multline}
  \frac{1}{h^2} \left(\frac{s_1}{s_2}\right)_{i-\hf,j} u_{i-1,j} + 
  \frac{1}{h^2} \left(\frac{s_1}{s_2}\right)_{i+\hf,j} u_{i+1,j} + 
  \frac{1}{h^2} \left(\frac{s_2}{s_1}\right)_{i,j-\hf} u_{i,j-1} + 
  \frac{1}{h^2} \left(\frac{s_2}{s_1}\right)_{i,j+\hf} u_{i,j+1} \\
  + \left(
    \frac{\omega^2}{(s_1s_2)_{i,j}\cdot c_{i,j}^2}
    - \left(\cdots
    \right) 
  \right)
  u_{i,j} = f_{i,j}
  \label{eq:helmd}
\end{multline}
with $u_{i',j'}$ equal to zero for $(i',j')$ that violates $1\le
i',j'\le n$. Here $(\cdots)$ stands for the sum of the four
coefficients that appear in the first line. We order both $u_{i,j}$
and $f_{i,j}$ row by row starting from the first row $j=1$ and define
the vectors
\begin{eqnarray*}
&&u = \left(u_{1,1}, u_{2,1},\ldots,u_{n,1}, \ldots,u_{1,n}, u_{2,n},\ldots,u_{n,n}\right)^t,\\
&&f = \left(f_{1,1}, f_{2,1},\ldots,f_{n,1}, \ldots,f_{1,n}, f_{2,n},\ldots,f_{n,n}\right)^t.
\end{eqnarray*}
Denote the discrete system of \eqref{eq:helmd} by $A u = f$. We
further introduce a block version of it by defining $\P_m$ to be set
of the indices in the $m$-th row
\[
\P_m = \{p_{1,m}, p_{2,m,} \ldots, p_{n,m}\}
\]
and introducing
\[
u_m = \left( u_{1,m}, u_{2,m},\ldots,u_{n,m} \right)^t \quad\text{and}\quad
f_m = \left( f_{1,m}, f_{2,m},\ldots,f_{n,m} \right)^t.
\]
Then
\[
u = (u_1^t, u_2^t, \ldots, u_n^t)^t,\quad
f = (f_1^t, f_2^t, \ldots, f_n^t)^t.
\]
Using these notations, the system $Au=f$ takes the following block
tridiagonal form
\[
\begin{pmatrix}
  A_{1,1} & A_{1,2} & & \\
  A_{2,1} & A_{2,2} & \ddots & \\
  & \ddots & \ddots & A_{n-1,n}\\
  & & A_{n,n-1} & A_{n,n}
\end{pmatrix}
\begin{pmatrix}
  u_1\\
  u_2\\
  \vdots\\
  u_n
\end{pmatrix}
=
\begin{pmatrix}
  f_1\\
  f_2\\
  \vdots\\
  f_n
\end{pmatrix}
\]
where $A_{m,m}$ are tridiagonal and $A_{m,m-1} = A_{m-1,m}^t$ are
diagonal matrices. 

The sweeping factorization of the matrix $A$ is essentially a block
$LDL^t$ factorization that eliminates the unknowns layer by layer,
starting from the absorbing layer near $x_2=0$. The result of this
process is a factorization
\begin{equation}
A = L_1 \cdots L_{n-1}
\begin{pmatrix}
  S_1 & & & \\
  & S_2 & & \\
  & & \ddots & \\
  & & & S_n\\
\end{pmatrix}
L_{n-1}^t \cdots L_1^t,
\label{eq:Afact}
\end{equation}
where $S_1=A_{1,1}$, $S_m = A_{m,m} - A_{m,m-1} S_{m-1}^{-1}
A_{m-1,m}$ for $m=2,\ldots, n$, and $L_k$ is given by
\[
L_k(\P_{k+1},\P_k) = A_{k+1,k} S_k^{-1},\quad
L_k(\P_i,\P_i) = I\;\;(1\le i \le n),\quad
\text{and zero otherwise}.
\]
This process is illustrated graphically in Figure \ref{fig:2Dgrid}
(right). Inverting this factorization for $A$ gives the following
formula for $u$:
\[
u = (L_1^t)^{-1}\cdots (L_{n-1}^t)^{-1} 
\begin{pmatrix}
  S_1^{-1} & & & \\
  & S_2^{-1} & & \\
  & & \ddots & \\
  & & & S_n^{-1} \\
\end{pmatrix}
L_{n-1}^{-1} \cdots L_1^{-1} f.
\]
Algorithmically, the construction of the sweeping factorization of $A$
can be summarized as follows by introducing $T_m = S_m^{-1}$.
\begin{algo}
  Construction of the sweeping factorization of $A$.
  \label{alg:setupext}
\end{algo}
\begin{algorithmic}[1]
  \STATE $S_1 = A_{1,1}$ and $T_1 = S_1^{-1}$.
  \FOR{$m=2,\ldots,n$}
  \STATE $S_m = A_{m,m} - A_{m,m-1} T_{m-1} A_{m-1,m}$ and $T_m = S_m^{-1}$.
  \ENDFOR
\end{algorithmic}
Since $S_m$ and $T_m$ are in general dense matrices of size $n\times
n$, the cost of the construction algorithm is of order $O(n^4) =
O(N^2)$. The computation of $u = A^{-1}f$ is carried out in the
following algorithm once the factorization is ready.
\begin{algo}
  Computation of $u=A^{-1}f$ using the sweeping factorization of $A$.
  \label{alg:solveext}
\end{algo}
\begin{algorithmic}[1]
  \FOR{$m=1,\ldots,n$}
  \STATE $u_m = f_m$
  \ENDFOR
  \FOR{$m=1,\ldots,n-1$}
  \STATE $u_{m+1} = u_{m+1} - A_{m+1,m} (T_m u_m)$
  \ENDFOR
  \FOR{$m=1,\ldots,n$}
  \STATE $u_m = T_m u_m$
  \ENDFOR
  \FOR{$m=n-1,\ldots,1$}
  \STATE $u_m = u_m - T_m (A_{m,m+1} u_{m+1})$
  \ENDFOR
\end{algorithmic}
Obviously the computations of $T_m u_m$ in the second and the third
loops only needs to be carried out once for each $m$. We prefer to
write the algorithm this way for the simplicity of presentation. The
cost of computing $u$ with Algorithm \ref{alg:solveext} is of order
$O(n^3) = O(N^{3/2})$, which is about $O(N^{1/2})$ times more
expensive compared to the multifrontal method. Therefore, these two
algorithms themselves are not very useful.

\subsection{Moving PML}

In Algorithms \ref{alg:setupext} and \ref{alg:solveext}, the dominant
cost is the construction and the application of the matrices $T_m$.
In \cite{EngquistYing:10a}, we emphasized the physical meaning of the
Schur complement matrices $T_m$ of the sweeping
factorization. Consider only the top-left $m\times m$ block of the
above factorization.
\begin{equation}
\begin{pmatrix}
  A_{1,1} & A_{1,2} & & \\
  A_{2,1} & A_{2,2} & \ddots & \\
  & \ddots & \ddots & A_{m-1,m}\\
  & & A_{m,m-1} & A_{m,m}
\end{pmatrix}
=
L_1 \cdots L_{m-1}
\begin{pmatrix}
  S_1 & & & \\
  & S_2 & & \\
  & & \ddots & \\
  & & & S_m\\
\end{pmatrix}
L_{m-1}^t \cdots L_1^t,
\label{eq:Afactm}
\end{equation}
where the $L_k$ matrices are redefined to their restriction to the
top-left $m \times m$ blocks. The matrix on the left is in fact the
discrete Helmholtz equation restricted to the half space below
$x_2=(m+1)h$ and with zero boundary condition on this line. Inverting
the factorization \eqref{eq:Afactm} gives
\[
\begin{pmatrix}
  A_{1,1} & A_{1,2} & & \\
  A_{2,1} & A_{2,2} & \ddots & \\
  & \ddots & \ddots & A_{m-1,m}\\
  & & A_{m,m-1} & A_{m,m}
\end{pmatrix}^{-1}
=
(L_1^t)^{-1} \cdots (L_{m-1}^t)^{-1}
\begin{pmatrix}
  S_1^{-1} & & & \\
  & S_2^{-1} & & \\
  & & \ddots & \\
  & & & S_m^{-1}\\
\end{pmatrix}
L_{m-1}^{-1} \cdots L_1^{-1}.
\]
The matrix on the left side is an approximation of the discrete
half-space Green's function of the Helmholtz operator with zero
boundary condition. On the right side, due to the definition of the
matrices $L_1,\ldots,L_{m-1}$, the $(m,m)$-th block of the product is
exactly equal to $S_m^{-1}$. Therefore,
\begin{center}
  $T_m = S_m^{-1}$ approximates the discrete half-space Green function
  of the
  Helmholtz operator \\
  with zero boundary at $x_2=(m+1)h$, restricted to the points on
  $x_2=mh$.
\end{center}

In the previous paper \cite{EngquistYing:10a}, $T_m$ is approximated
using the hierarchical matrix framework. Due to the fact that the 3D
Green's function, restricted to a plane, propagates oscillations in
all directions, the theoretical justification of that method is
lacking in 3D. Here, we try to approximate the matrix $T_m$ in a
different way. 


As an operator, $T_m: g_m \rightarrow v_m$ maps an external force
$g_m$ loaded only on the $m$-th layer to the solution $v_m$ restricted
to the same layer. Though it is a map between quantities only defined
on the $m$-th layer, the computation domain includes all first $m$
layers with the PML padded near $x_2=0$. However, since the force
$g_m$ is only loaded on the $m$-th layer, there is no reason to keep
the PML layer near $x_2=0$ if one can be satisfied with an
approximation. 
\begin{center}
  The central idea is to push the PML from $x_2=0$ right next to
  $x_2=mh$.
\end{center}
To make this precise, let us assume that the width $\eta$ of the PML
is an integer multiple of $h$ and let $b=\eta/h$ be the number of grid
points in PML layer in the transversal direction. Define
\[
s^m_2(x_2) = \left( 1+i\frac{\sigma_2(x_2-(m-b)h)}{\omega} \right)^{-1}
\]
and introduce an auxiliary problem on the domain $D_m = [0,1]\times
[(m-b)h,(m+1)h]$:
\begin{eqnarray}
\left( (s_1\p_1)(s_1\p_1) + (s^m_2\p_2)(s^m_2\p_2) + \frac{\omega^2}{c^2(x)} \right) u = f && x\in D_m, \label{eq:helma}\\
u = 0 && x \in \p D_m . \nonumber
\end{eqnarray}
This equation is discretized with the subgrid
\[
\G_m = \{p_{i,j}, 1\le i \le n, m-b+1 \le j \le m \}
\]
of the original grid $\P$ and the resulting $bn \times bn$ discrete
Helmholtz operator is denoted by $H_m$. Following the main idea
mentioned above, the operator $\wt{T}_m: g_m \rightarrow v_m$ defined
through $H_m$ by
\[
\begin{pmatrix}
  *\\ \vdots \\ * \\ v_m
\end{pmatrix}
\approx 
H_m^{-1}
\begin{pmatrix}
  0\\ \vdots \\ 0 \\ g_m
\end{pmatrix}
\]
is an approximation to the matrix $T_m$. Notice that applying
$\wt{T}_m$ to an arbitrary vector $g_m$ involves solving a linear
system of matrix $H_m$, which comes from the {\em local} 5-point
stencil on the narrow grid $\G_m$ that contains only $b$ layers. Let
us introduce a new ordering for $\G_m$
\[
p_{1,m-b+1},p_{1,m-b+2},\ldots,p_{1,m}\ldots
p_{n,m-b+1},p_{n,m-b+2},\ldots,p_{n,m}
\]
that iterates through the $x_2$ direction and denote the permutation
matrix induced from this new ordering by $P_m$. Now the matrix $P_m
H_m P_m^t$ is a banded matrix with only $b-1$ lower diagonals and
$b-1$ upper diagonals. It is well known that the LU factorization $L_m
U_m = P_m H_m P_m^t$ can be constructed efficiently. As a result, the
application of $\wt{T}_m$ can be carried out rapidly.

We call this approach the {\em moving PML} method, since these new
PMLs do not exist in the original problem as they are only introduced
in order to approximate $T_m$ efficiently. In the above discussion,
the moving PML is pushed right next to $x_2=mh$. However, in general
we can place the moving PML at a location that is a few layers away
from $x_2=mh$. The potential advantage of keeping a few extra layers
as a buffer is that the resulting approximation $\wt{T}_m$ is more
accurate. On the other hand, since there are more layers in the
subgrid $\G_m$ for each $m$, the computational cost grows
accordingly. In our numerical tests, we observe that extra buffer
layers provide little improvement on the approximation accuracy and
hence the moving PML is indeed pushed right next to $x_2=mh$.

The application of a PML right next to the layer to be eliminated
corresponds to a PML or absorbing boundary condition next to a
Dirichlet boundary condition. This has been used as an asymptotic
technique for high frequency scattering under the name of on-surface
radiation boundary condition (OSRBC)
\cite{AtleEngquist:07,KriegsmannTaflove:87}. The OSRBC is an
approximation that is more accurate than physical optics but, of
course, not as accurate as a full boundary integral formulation.

\subsection{Approximate inversion and preconditioner}

Let us incorporate the moving PML technique into Algorithms
\ref{alg:setupext} and \ref{alg:solveext}. The computation at the
first $(b+1)$ layers needs to be handled differently, since it does
not make sense to introduce moving PML for these initial layers. Let
us call the first $b$ layers the {\em front} part and define
\[
u_F = (u_1^t,\ldots,u_{b}^t)^t \quad \text{and} \quad
f_F = (f_1^t,\ldots,f_{b}^t)^t.
\]
Then we can rewrite $Au=f$ as
\[
\begin{pmatrix}
  A_{F,F} & A_{F,b+1} & & \\
  A_{b+1,F} & A_{b+1,b+1} & \ddots & \\
  & \ddots & \ddots & A_{n-1,n}\\
  & & A_{n,n-1} & A_{n,n}
\end{pmatrix}
\begin{pmatrix}
  u_F\\
  u_{b+1}\\
  \vdots\\
  u_n
\end{pmatrix}
=
\begin{pmatrix}
  f_F\\
  f_{b+1}\\
  \vdots\\
  f_n
\end{pmatrix}.
\]
The construction of the approximate sweeping factorization of $A$
takes the following steps. Notice that since $T_m$ are approximated
directly there is no need to compute $S_m$ anymore.
\begin{algo}
  Construction of the approximate sweeping factorization of $A$ with
  moving PML.
  \label{alg:setup}
\end{algo}
\begin{algorithmic}[1]
  \STATE Let $\G_F$ be the subgrid of the first $b$ layers, $H_F =
  A_{F,F}$, and $P_F$ be the permutation induce by the new ordering
  ($x_2$ first) of $\G_F$. Construct the LU factorization $L_F U_F =
  P_F H_F P_F^t$. This factorization implicitly defines $\wt{T}_F:
  \C^{bn} \rightarrow \C^{bn}$.

  \FOR{$m=b+1,\ldots,n$}

  \STATE Let $\G_m = \{p_{i,j}, 1\le i \le n, m-b+1 \le j \le m \}$,
  $H_m$ be the discrete system of \eqref{eq:helma} on $\G_m$, and
  $P_m$ be the permutation induced by the new ordering of
  $\G_m$. Construct the LU factorization $L_m U_m = P_m H_m
  P_m^t$. This factorization implicitly defines $\wt{T}_m: \C^n
  \rightarrow \C^n$.
  
  \ENDFOR
\end{algorithmic}
The cost of Algorithm \ref{alg:setup} is $O(b^3 n^2) = O(b^3 N)$. The
computation of $u \approx A^{-1} f$ using the constructed sweeping
factorization is summarized in the following algorithm
\begin{algo}
  Computation of $u \approx A^{-1}f$ using the sweeping factorization
  of $A$ with moving PML.
  \label{alg:solve}
\end{algo}
\begin{algorithmic}[1]
  \STATE $u_F = f_F$ and $u_m = f_m$ for $m=b+1,\ldots,n$.

  \STATE $u_{b+1} = u_{b+1} - A_{b+1,F} (\wt{T}_F u_F)$. $\wt{T}_F
  u_F$ is computed as $P_F^t U_F^{-1} L_F^{-1} P_F u_F$.
  
  \FOR{$m=b+1,\ldots,n-1$}
  
  \STATE $u_{m+1} = u_{m+1} - A_{m+1,m} (\wt{T}_m u_m)$. The
  application of $\wt{T}_m u_m$ is done by forming the vector
  $(0,\ldots,0,u_m^t)^t$, applying $P_m^t U_m^{-1} L_m^{-1} P_m$ to
  it, and extracting the value on the last layer.
  
  \ENDFOR
  
  \STATE $u_F = \wt{T}_F u_F$. See the previous steps for the
  application of $\wt{T}_F$.
  
  \FOR{$m=b+1,\ldots,n$}
  
  \STATE $u_m = \wt{T}_m u_m$. See the previous steps for the
  application of $\wt{T}_m$.
  
  \ENDFOR
  
  \FOR{$m=n-1,\ldots,b+1$}
  
  \STATE $u_m = u_m - \wt{T}_m (A_{m,m+1} u_{m+1})$. See the previous
  steps for the application of $\wt{T}_m$.
  
  \ENDFOR

  \STATE $u_F = u_F - \wt{T}_F (A_{F,b+1} u_{b+1})$. See the previous
  steps for the application of $\wt{T}_F$.
  
\end{algorithmic}
The cost of Algorithm \ref{alg:solve} is $O(b^2 n^2) = O(b^2
N)$. Since $b$ is a fixed constant, the cost is essentially linear.
Algorithm \ref{alg:solve} defines an operator
\[
M: 
f = (f_F^t, f_{b+1}^t, \ldots, f_n^t)^t \rightarrow
u = (u_F^t, u_{b+1}^t, \ldots, u_n^t)^t,
\]
which is an approximate inverse of the discrete Helmholtz operator
$A$. Due to the indefiniteness of $A$, this approximate inverse might
suffer from instability. In practice, instead of generating the
sweeping factorization of the original matrix $A$, we choose to
generate the factorization for the matrix $A_\alpha$ associated with
the modified Helmholtz equation
\begin{equation}
\Lapl u(x) + \frac{(\omega+i\alpha)^2}{c^2(x)} u(x) = f(x),
\label{eq:helmalpha}
\end{equation}
where $\alpha$ is an $O(1)$ positive constant. We denote by $M_\alpha:
f \rightarrow u$ the operator defined by Algorithm \ref{alg:solve}
with this modified equation. We would like to emphasize that
\eqref{eq:helmalpha} is very different from the equation used in the
shifted Laplacian approach (for example
\cite{ErlanggaOosterleeVuik:06,LairdGiles:02}): in the shifted
Laplacian formulation the imaginary part of the operator is
$O(\omega)$ while here the imaginary part is $O(1)$.

Since $\alpha$ is small, $A_\alpha$ is close to $A$. Therefore, we
propose to solve the preconditioner system
\[
M_\alpha A u = M_\alpha f
\]
using the GMRES solver \cite{Saad:03,SaadSchultz:86}. As the cost of
applying $M_\alpha$ to any vector is $O(n^2)=O(N)$, the total cost of
the iterative solver scales like $O(N_I N)$, where $N_I$ is the number
of iterations. As the numerical results in Section \ref{sec:2Dnum}
demonstrate, $N_I$ depends at most logarithmically on $N$, thus
resulting a solver with almost linear complexity.

\begin{figure}[h!]
  \begin{center}
    \includegraphics{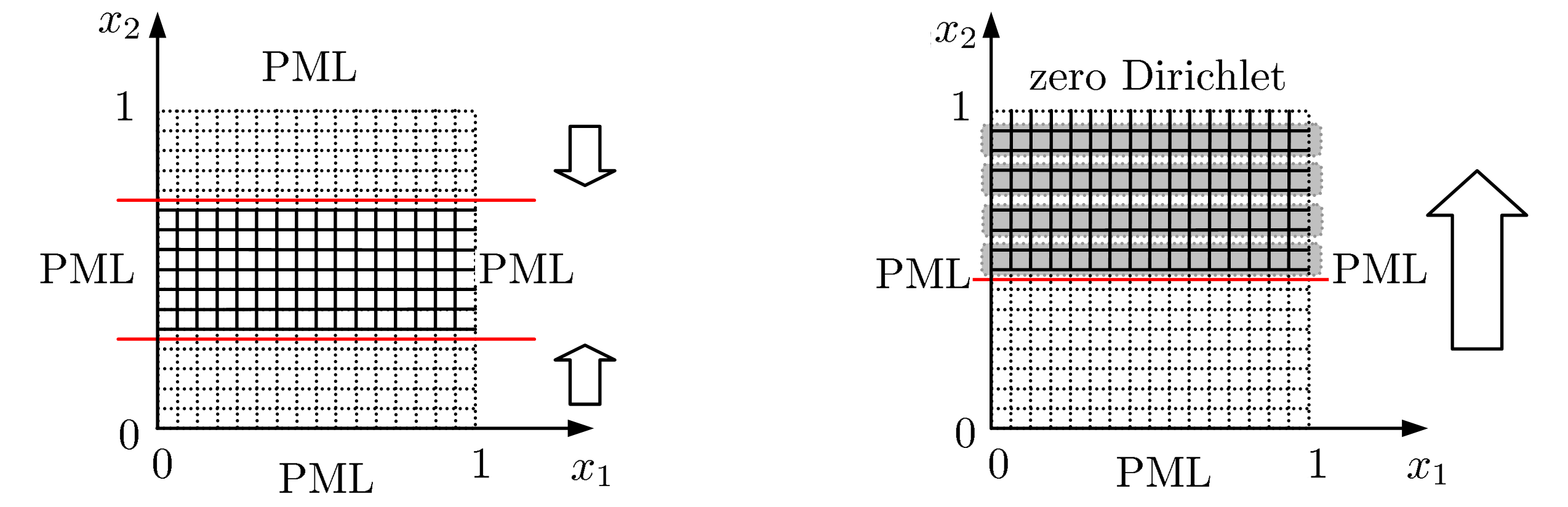}
  \end{center}
  \caption{Different sweeping patterns. Left: For problems with PML at
    both $x_2=0$ and $x_2=1$, the algorithm sweeps from both ends
    towards the center. Right: Instead of one layer, multiple layers
    of unknowns can be eliminated within each iteration of the
    algorithm.}
  \label{fig:2Dvar}
\end{figure}

The problem considered so far has zero Dirichlet boundary condition on
$x_2=1$. A common situation is to impose PML at all sides. In this
case, the algorithms need a slight modification. Instead of sweeping
upward from $x_2=0$, the algorithm sweeps with two fronts, one from
$x_2=0$ upward and the other from $x_2=1$ downward (see Figure
\ref{fig:2Dvar} (left)). Similar to $u_F$ and $f_F$ near $x_2=0$, we
introduce
\[
u_E = (u_{n-b+1}^t,\ldots,u_n^t)^t \quad\text{and}\quad
f_E = (f_{n-b+1}^t,\ldots,f_n^t)^t
\]
and write $Au=f$ in the following block form
\[
\begin{pmatrix}
  A_{F,F} & A_{F,b+1} & & &\\
  A_{b+1,F} & A_{b+1,b+1} & \ddots & &\\
  & \ddots & \ddots & \ddots & \\
  & & \ddots & A_{n-b,n-b} & A_{n-b,E}\\
  & & & A_{E,n-b} & A_{E,E}
\end{pmatrix}
\begin{pmatrix}
  u_F\\
  u_{b+1}\\
  \vdots\\
  u_{n-b}\\
  u_E
\end{pmatrix}
=
\begin{pmatrix}
  f_F\\
  f_{b+1}\\
  \vdots\\
  f_{n-b}\\
  f_E
\end{pmatrix}.
\]
The upward sweep goes through $m=F,b+1,\ldots,(n-1)/2$, and the
downward sweep visits $m=E,n-b,\ldots,(n+3)/2$. Finally, the
algorithm visits the middle layer $m=(n+1)/2$ with moving PMLs on
both sides.

Algorithm \ref{alg:setup} eliminates one layer of unknowns within each
iteration. We can also instead eliminate several layers of unknowns
together within each iteration (see Figure \ref{fig:2Dvar}
(right)). The resulting algorithm spends more computational time
within each elimination step, since the discrete system $H_m$ contains
more layers in the $x_2$ dimension. On the other hand, the number of
elimination steps goes down by a factor equal to the number of layers
processed within each elimination step. In practice, the actual number
$d$ of layers processed within each step depends on the width of the
moving PML and is chosen to minimize the overall computation time and
storage.

\section{Numerical Results in 2D}
\label{sec:2Dnum}

In this section, we present several numerical results to illustrate
the properties of the sweeping preconditioner described in Section
\ref{sec:2Dpre}. The algorithms are implemented in Matlab and all
tests are performed on a computer with a 2.6GHz CPU. We use GMRES as
the iterative solver with relative residue tolerance equal to
$10^{-3}$.

\subsection{PML}

\begin{figure}[h!]
  \begin{center}
    \begin{tabular}{ccc}
      \includegraphics[height=1.6in]{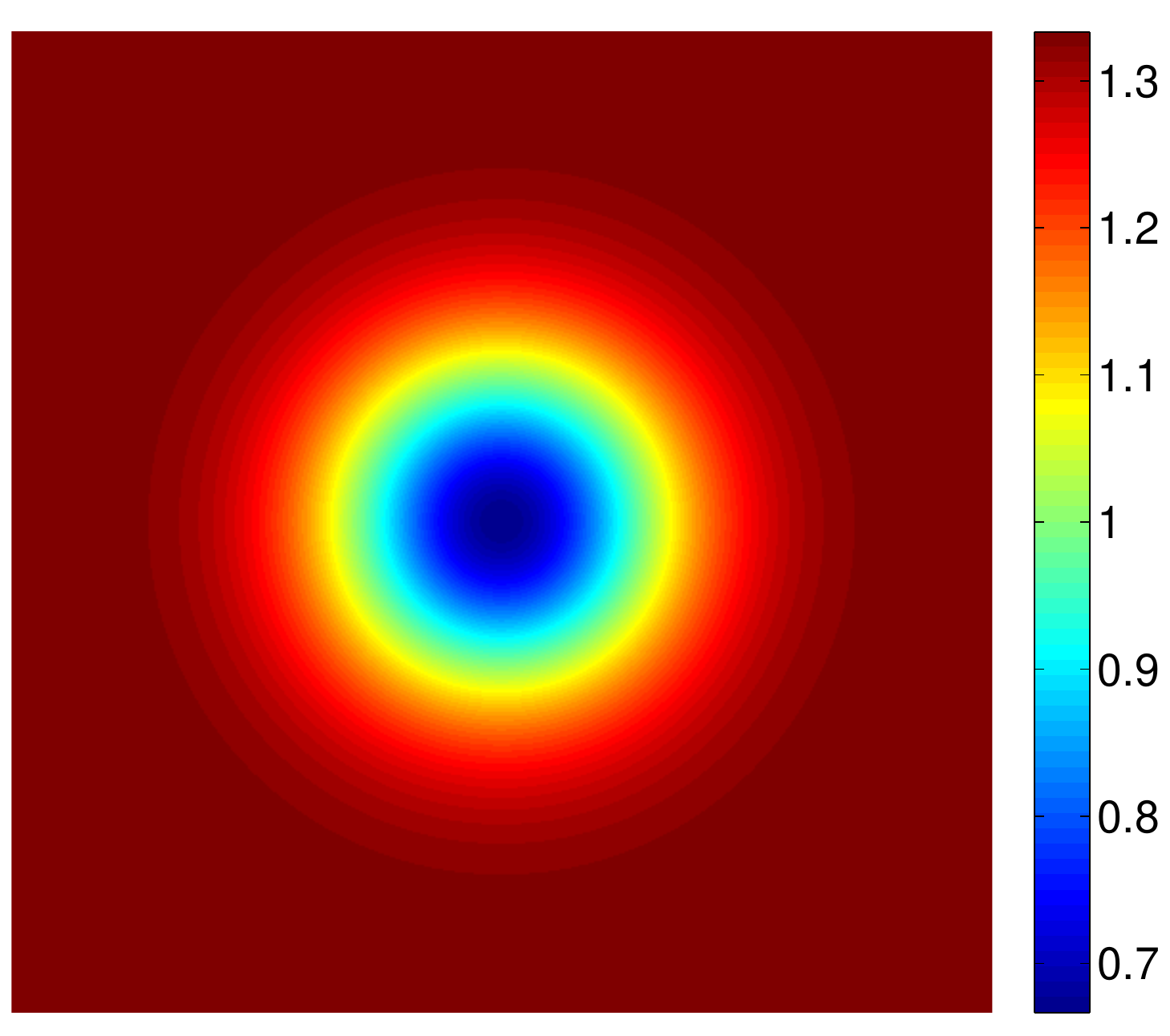}&\includegraphics[height=1.6in]{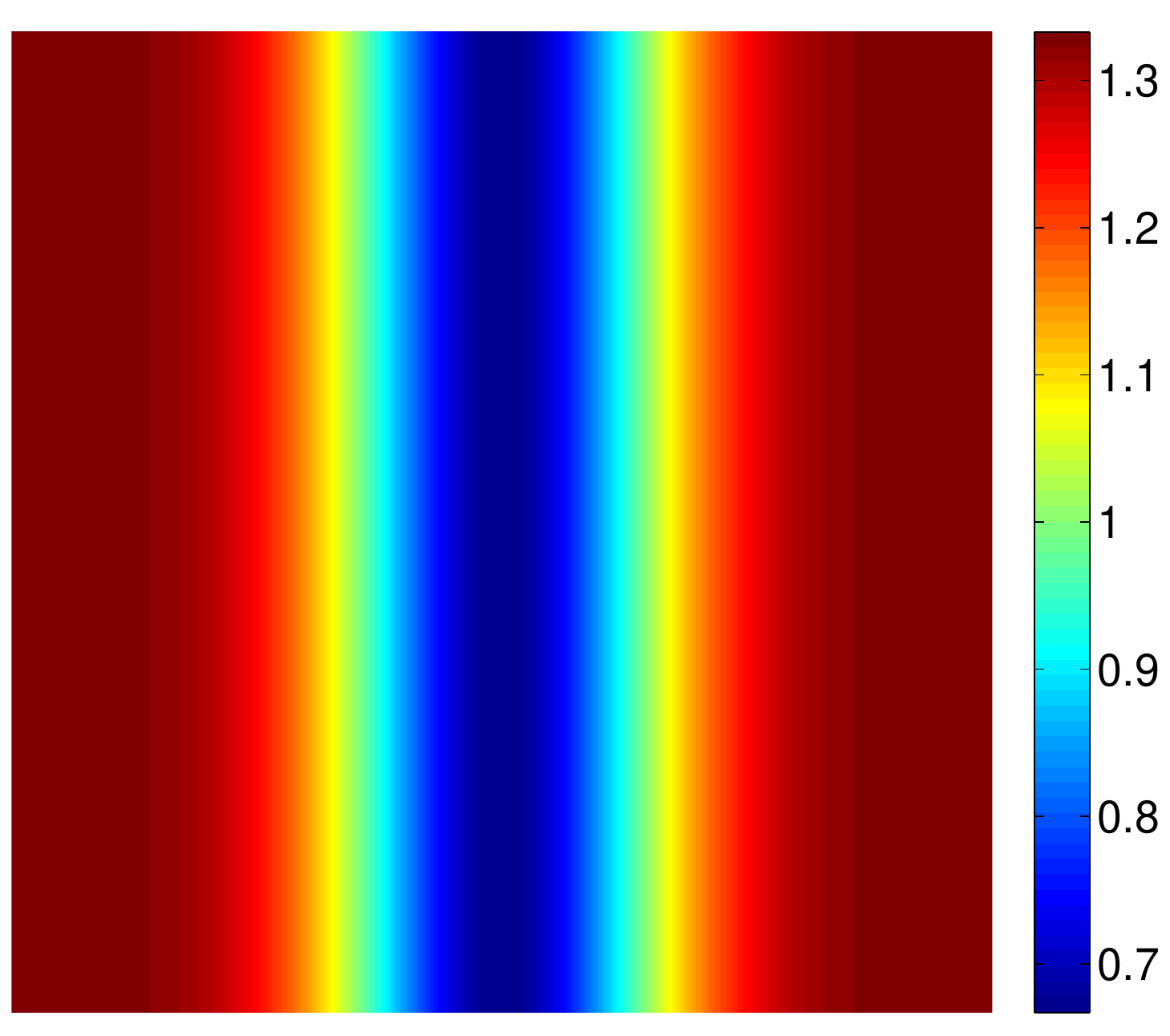}&\includegraphics[height=1.6in]{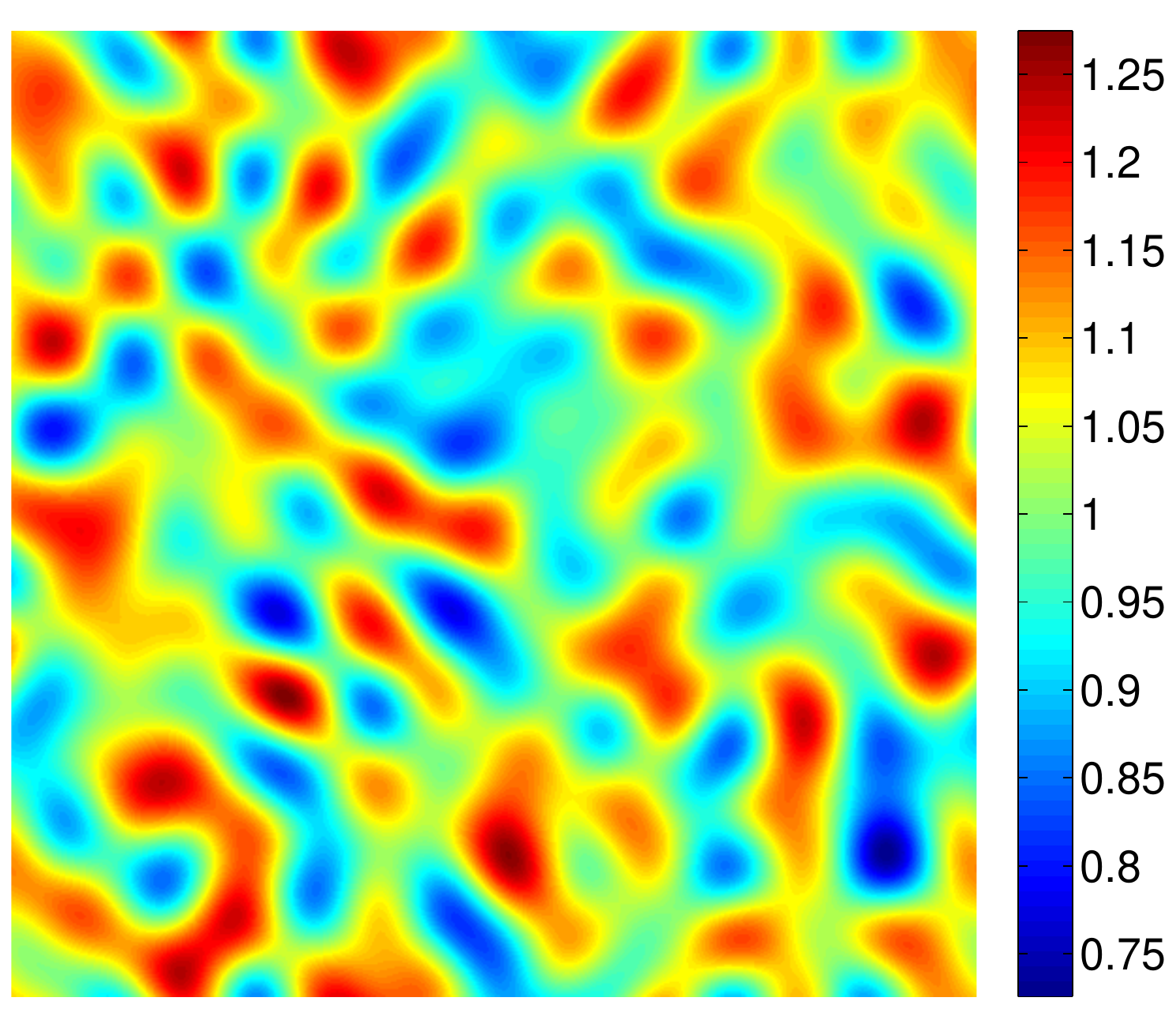}\\
      (a) & (b) & (c)
    \end{tabular}
  \end{center}
  \caption{Test velocity fields.}
  \label{fig:2Dnumspeed}
\end{figure}

The examples in this seciton have the PML boundary condition specified
at all sides. We consider three velocity fields in the domain
$D=(0,1)^2$:
\begin{enumerate}
\item The first velocity field corresponds to a smooth converging lens
  with a Gaussian profile at the center of the domain (see Figure
  \ref{fig:2Dnumspeed}(a)).
\item The second velocity field is a vertical waveguide with Gaussian cross
  section (see Figure \ref{fig:2Dnumspeed}(b)).
\item The third velocity field has a random velocity field (see Figure
  \ref{fig:2Dnumspeed}(c)).
\end{enumerate}

For each velocity field, we test with two external forces $f(x)$.
\begin{enumerate}
\item The first external force $f(x)$ is a narrow Gaussian point
  source located at $(x_1,x_2) = (0.5, 0.125)$. The response of this
  forcing term generates circular waves propagating at all
  directions. Due to the variations of the velocity field, the
  circular waves are going to bend, form caustics, and intersect.
\item The second external force $f(x)$ is a Gaussian wave packet whose
  wavelength is comparable to the typical wavelength of the
  domain. This packet centers at $(x_1,x_2) = (0.125, 0.125)$ and
  points to the $(1,1)$ direction. The response of this forcing term
  generates a Gaussian beam initially pointing towards the $(1,1)$
  direction. Due to the variations of the velocity field, this
  Gaussian beam bends and scatters.
\end{enumerate}

Firstly, we study how the sweeping preconditioner behaves when
$\omega$ varies. For each velocity field, we perform tests for
$\frac{\omega}{2\pi} = 16,32,\ldots,256$. In these tests, we
discretize each wavelength with $q=8$ points and $n=127, 255, \ldots,
2047$. The $\alpha$ value of the modified system is set to be equal to
2. The width of the moving PML is equal to $12h$ (i.e. $b=12$) and the
number $d$ of layers processed within each iteration of Algorithms
\ref{alg:setup} and \ref{alg:solve} is equal to 12. The sweeping
pattern indicated in Figure \ref{fig:2Dvar} (left) is used in these
tests.

\begin{table}[h!]
  \begin{center}
    \includegraphics[height=2.3in]{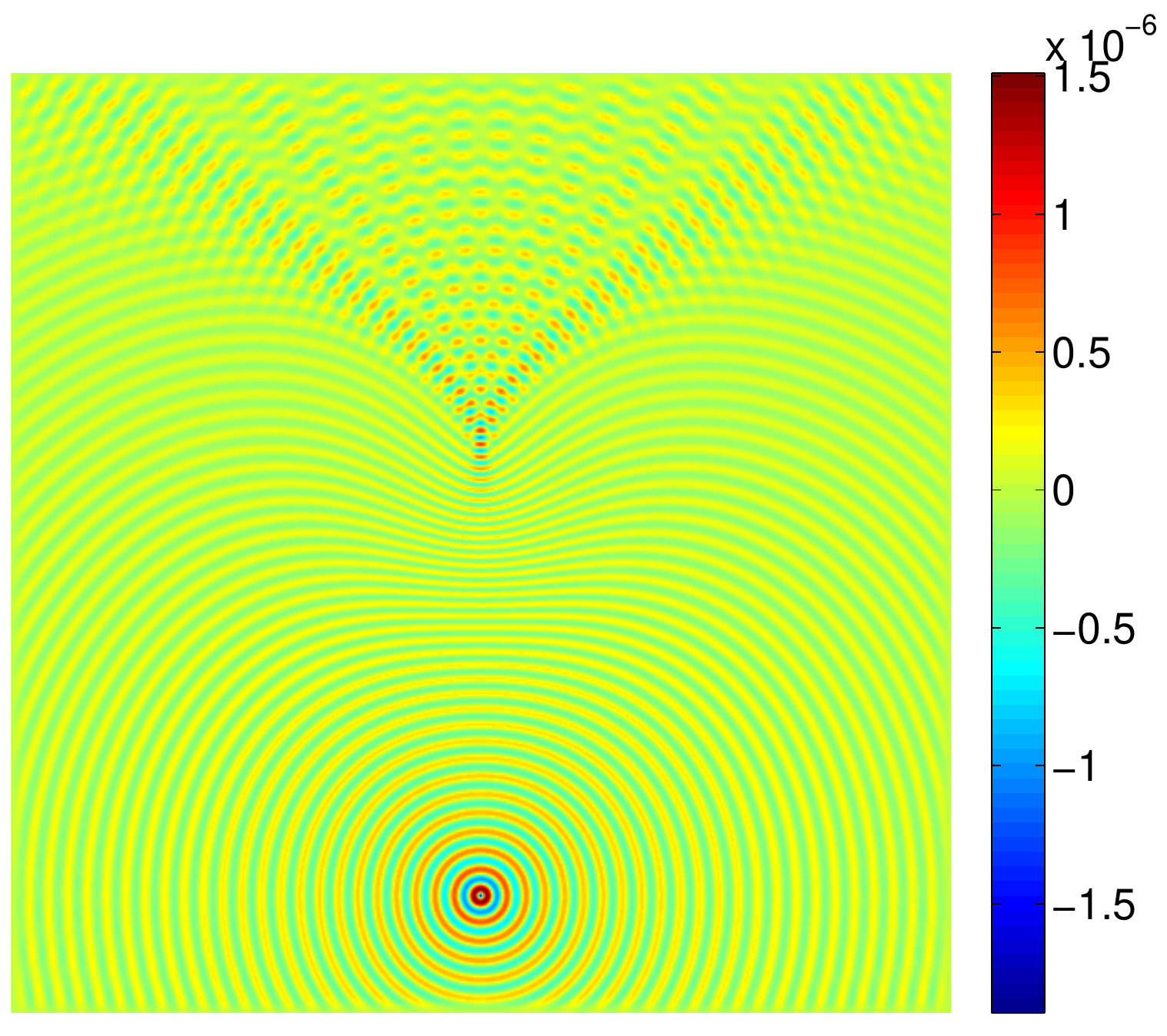}
    \includegraphics[height=2.3in]{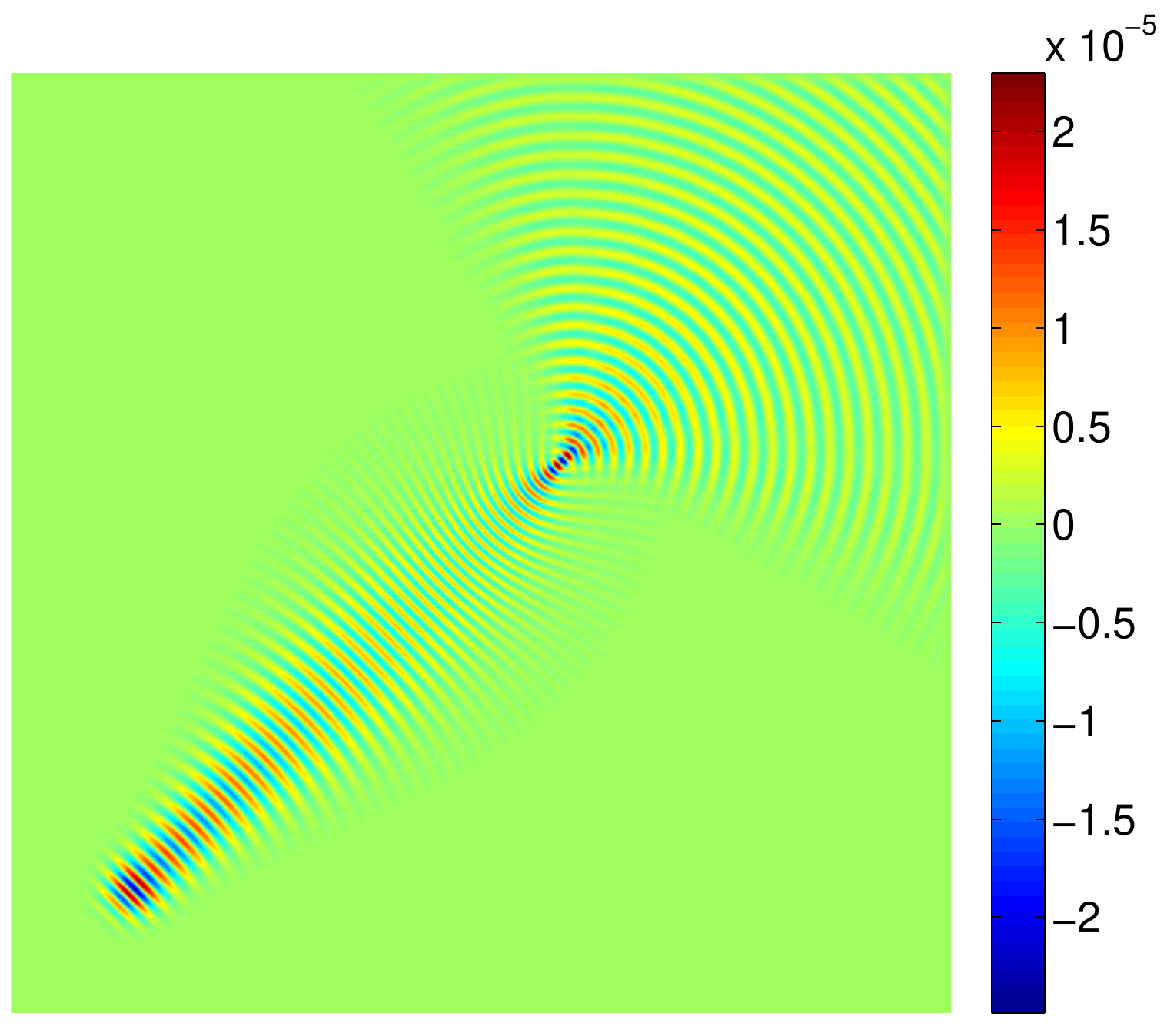}\\
    \begin{tabular}{|cccc|cc|cc|}
      \hline
      \multicolumn{4}{|c|}{} & \multicolumn{2}{|c|}{Test 1} & \multicolumn{2}{|c|}{Test 2}\\
      \hline
      $\omega/(2\pi)$ & $q$ & $N=n^2$ & $T_{\text{setup}}$ & $N_{\text{iter}}$ & $T_{\text{solve}}$ & $N_{\text{iter}}$ & $T_{\text{solve}}$\\
      \hline
      16  & 8 & $127^2$  & 2.86e-01 & 14 & 4.73e-01 & 15 & 3.81e-01\\
      32  & 8 & $255^2$  & 8.95e-01 & 15 & 1.59e+00 & 15 & 1.57e+00\\
      64  & 8 & $511^2$  & 3.78e+00 & 15 & 7.14e+00 & 15 & 7.12e+00\\
      128 & 8 & $1023^2$ & 1.61e+01 & 15 & 2.90e+01 & 13 & 2.54e+01\\
      256 & 8 & $2047^2$ & 6.85e+01 & 16 & 1.44e+02 & 11 & 9.42e+01\\
      \hline
    \end{tabular}
  \end{center}
  \caption{Results of velocity field 1 with varying $\omega$. Top: Solutions for two external forces with $\omega/(2\pi)=64$.
    Bottom: Results for different $\omega$.  }
  \label{tbl:2DPML1}
\end{table}

\begin{table}[h!]
  \begin{center}
    \includegraphics[height=2.3in]{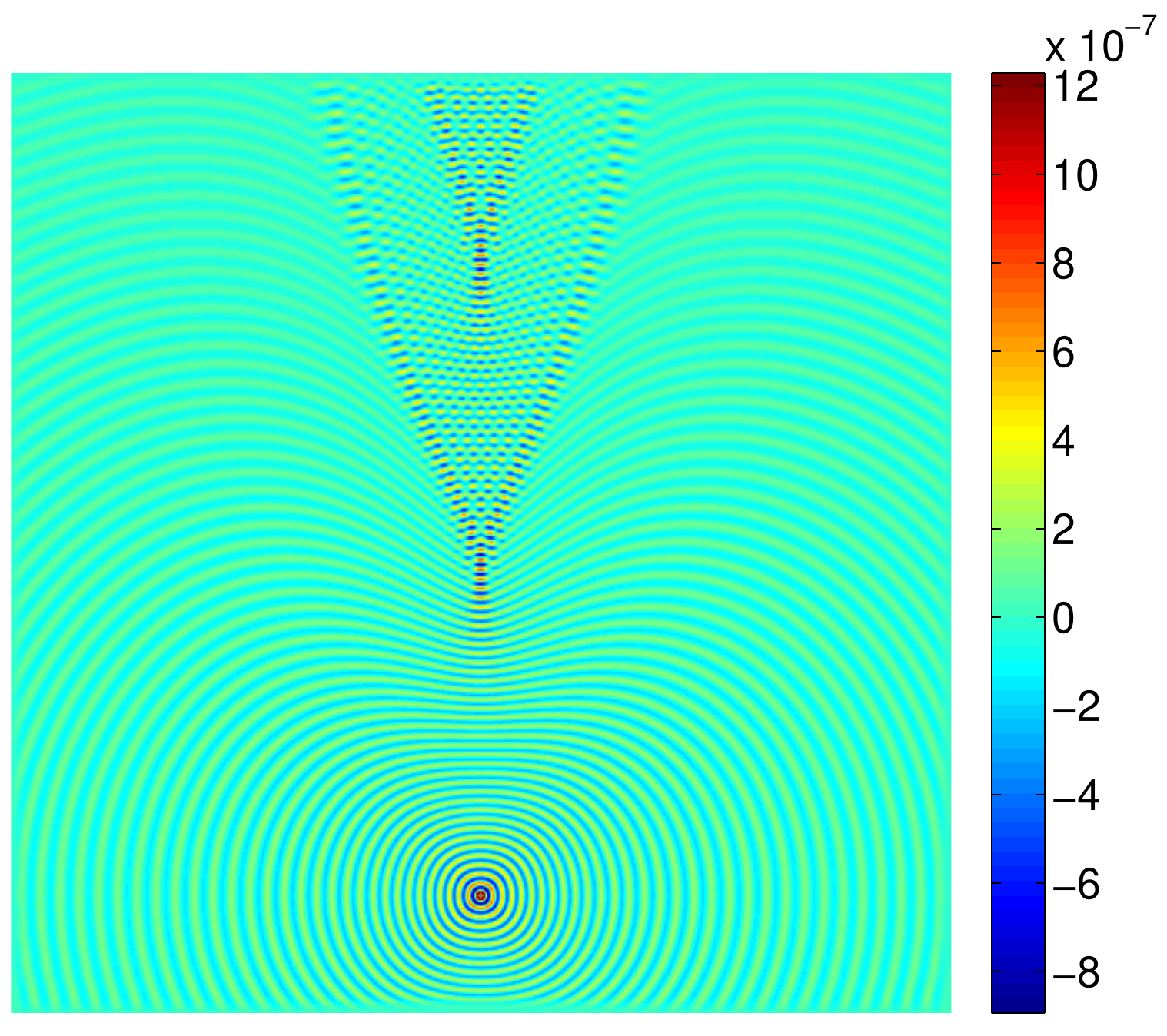}
    \includegraphics[height=2.3in]{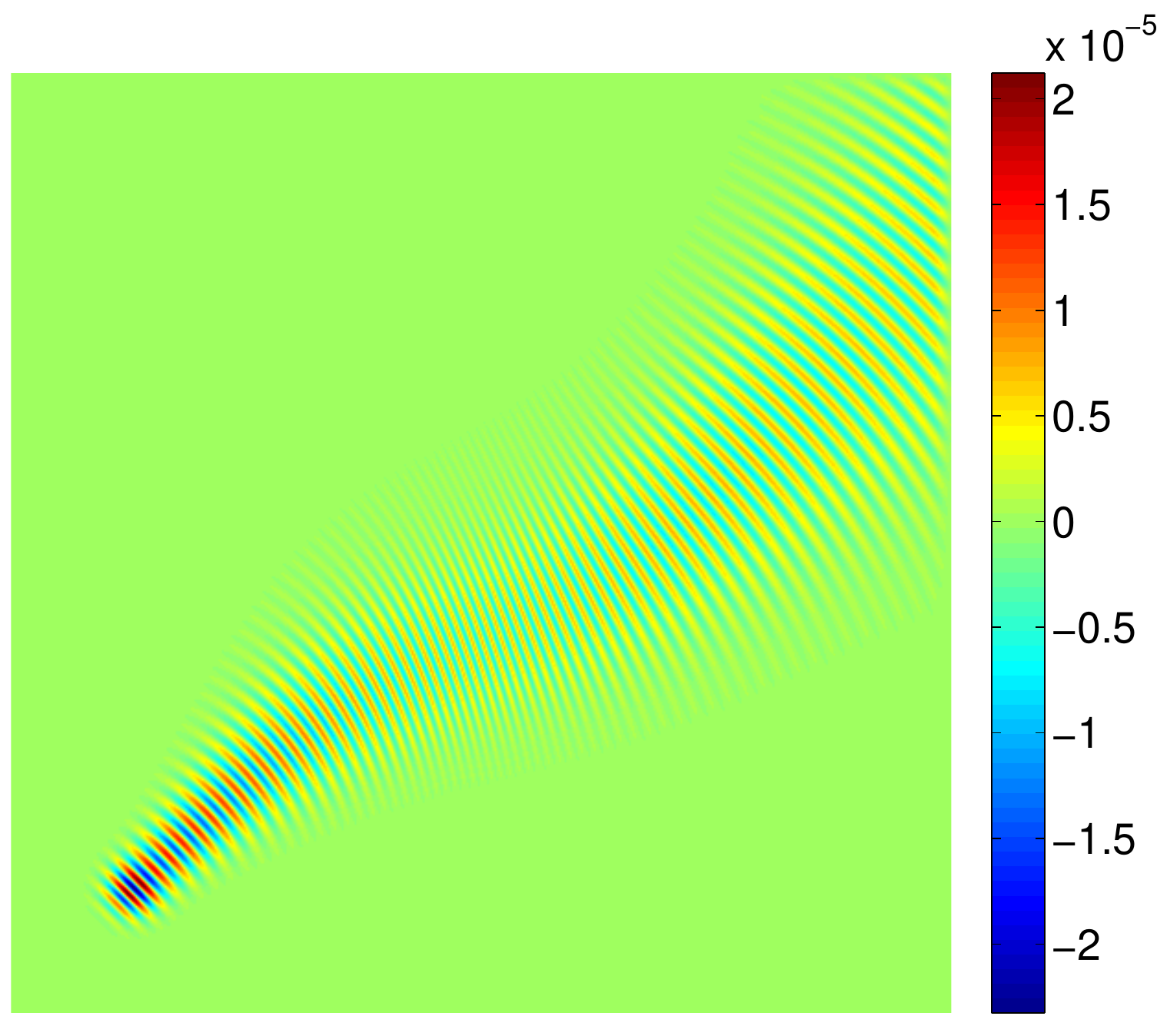}\\
    \begin{tabular}{|cccc|cc|cc|}
      \hline
      \multicolumn{4}{|c|}{} & \multicolumn{2}{|c|}{Test 1} & \multicolumn{2}{|c|}{Test 2}\\
      \hline
      $\omega/(2\pi)$ & $q$ & $N=n^2$ & $T_{\text{setup}}$ & $N_{\text{iter}}$ & $T_{\text{solve}}$ & $N_{\text{iter}}$ & $T_{\text{solve}}$\\
      \hline
      16  & 8 & $127^2$  & 3.03e-01 & 18 & 5.13e-01 & 16 & 3.94e-01\\
      32  & 8 & $255^2$  & 9.30e-01 & 19 & 2.58e+00 & 16 & 2.17e+00\\
      64  & 8 & $511^2$  & 3.76e+00 & 19 & 1.02e+01 & 15 & 7.61e+00\\
      128 & 8 & $1023^2$ & 1.61e+01 & 19 & 4.18e+01 & 13 & 2.75e+01\\
      256 & 8 & $2047^2$ & 6.78e+01 & 19 & 1.86e+02 & 12 & 1.10e+02\\

      \hline
    \end{tabular}
  \end{center}
  \caption{Results of velocity field 2 with varying $\omega$. Top: Solutions for two external forces with $\omega/(2\pi)=64$.
    Bottom: Results for different $\omega$.  }
  \label{tbl:2DPML2}
\end{table}

\begin{table}[h!]
  \begin{center}
    \includegraphics[height=2.3in]{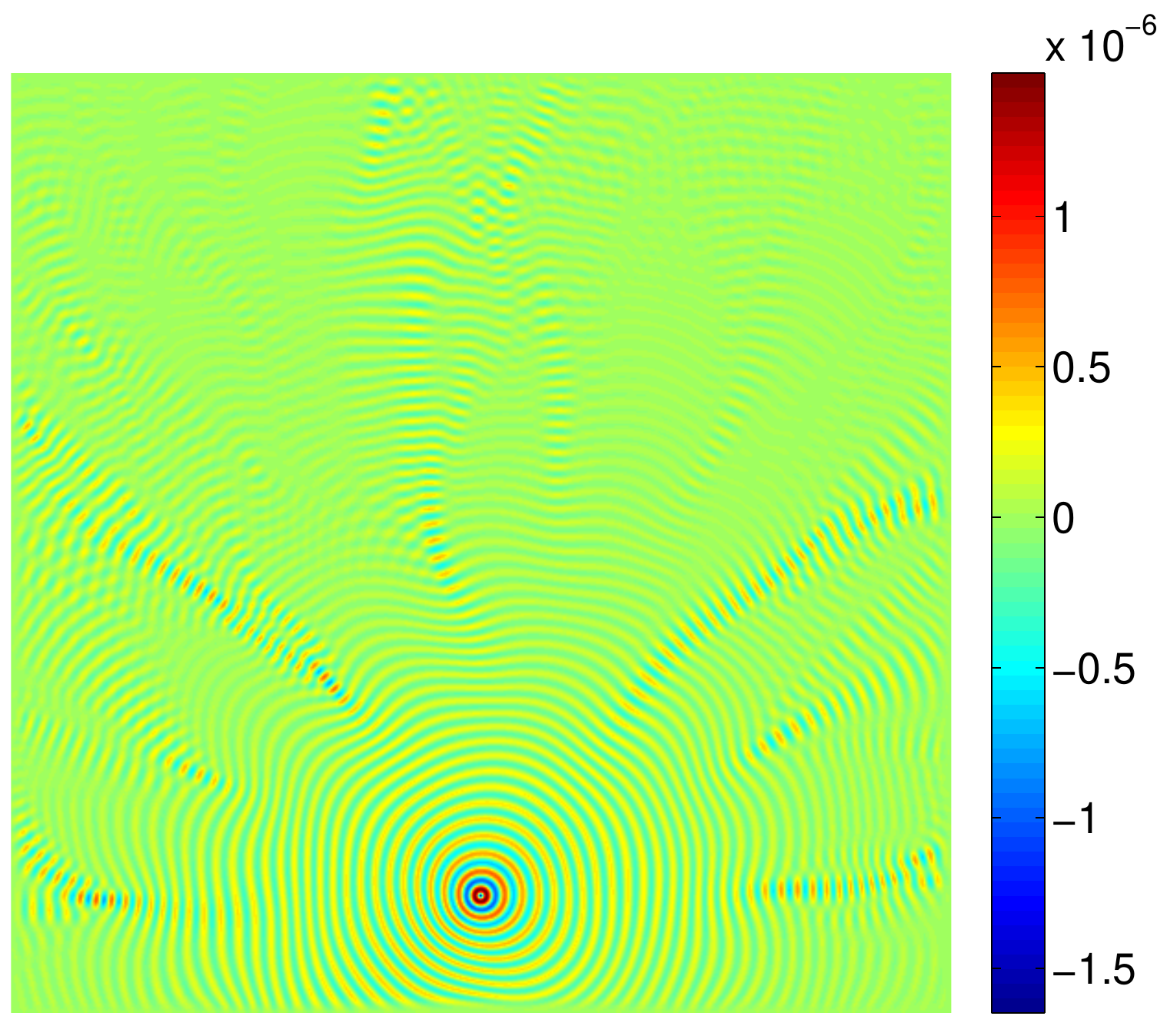}
    \includegraphics[height=2.3in]{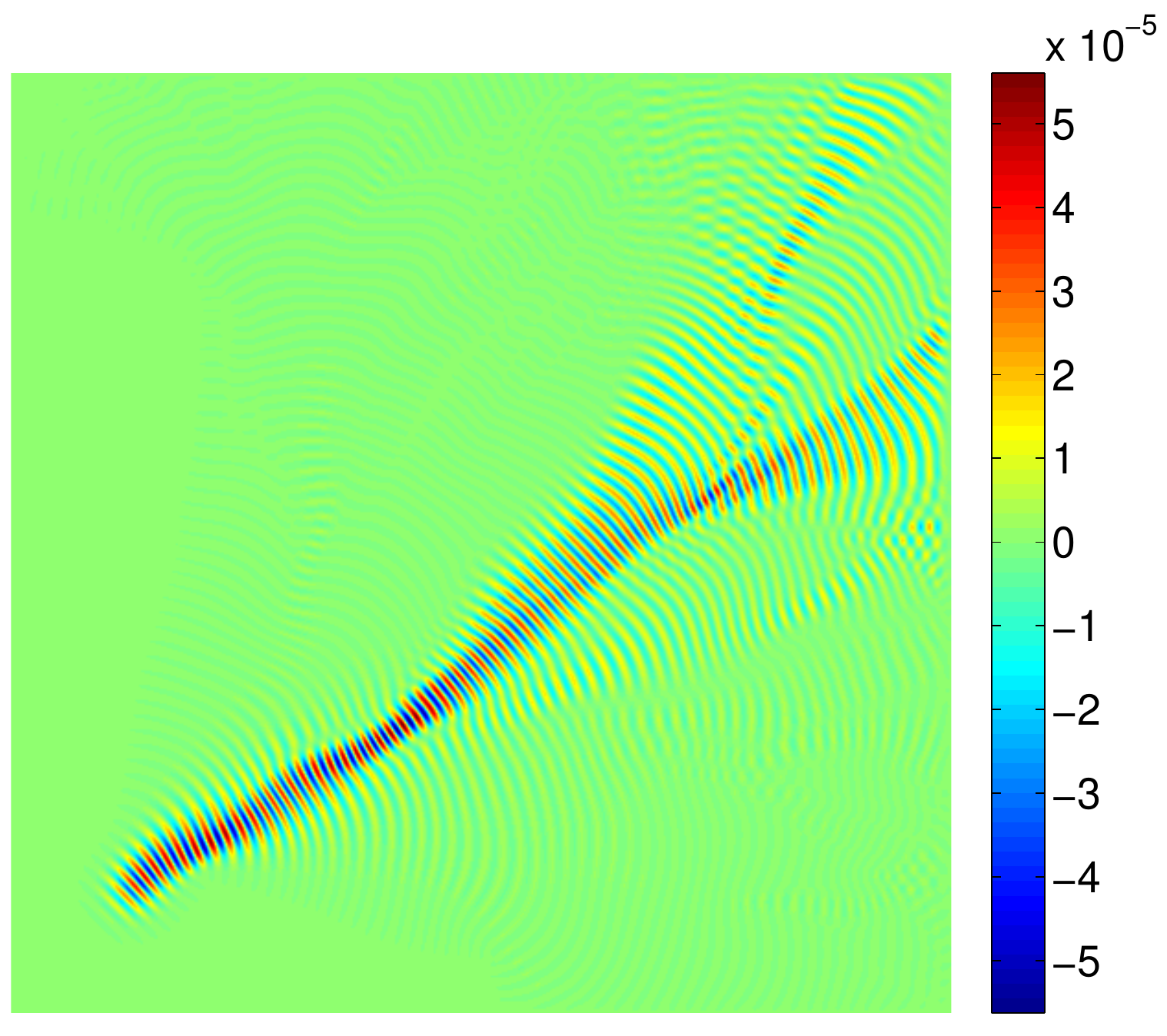}\\
    \begin{tabular}{|cccc|cc|cc|}
      \hline
      \multicolumn{4}{|c|}{} & \multicolumn{2}{|c|}{Test 1} & \multicolumn{2}{|c|}{Test 2}\\
      \hline
      $\omega/(2\pi)$ & $q$ & $N=n^2$ & $T_{\text{setup}}$ & $N_{\text{iter}}$ & $T_{\text{solve}}$ & $N_{\text{iter}}$ & $T_{\text{solve}}$\\
      \hline
      16  & 8 & $127^2$  & 3.56e-01 & 18 & 4.91e-01 & 19 & 5.03e-01\\
      32  & 8 & $255^2$  & 9.31e-01 & 18 & 2.42e+00 & 19 & 2.61e+00\\
      64  & 8 & $511^2$  & 3.76e+00 & 17 & 8.66e+00 & 23 & 1.24e+01\\
      128 & 8 & $1023^2$ & 1.60e+01 & 19 & 3.90e+01 & 22 & 4.80e+01\\
      256 & 8 & $2047^2$ & 6.82e+01 & 17 & 1.54e+02 & 17 & 1.48e+02\\
      \hline
    \end{tabular}
  \end{center}
  \caption{Results of velocity field 3 with varying $\omega$. Top: Solutions for two external forces with $\omega/(2\pi)=64$.
    Bottom: Results for different $\omega$.  }
  \label{tbl:2DPML3}
\end{table}

The results of the first velocity field are summarized in Table
\ref{tbl:2DPML1}. $T_{\text{setup}}$ denotes the time used to
construct the preconditioner in seconds. For each external force,
$N_{\text{iter}}$ is the number of iterations of the preconditioned
GMRES iteration and $T_{\text{solve}}$ is the overall solution
time. When $\omega$ and $n$ double, $N$ increases by a factor of 4 and
the setup cost in Table \ref{tbl:2DPML1} increases roughly by a factor
of $4$ as well, which is consistent with the $O(N)$ complexity of
Algorithm \ref{alg:setup}. At the same time, the number of iterations
is essentially independent of $n$. As a result, the overall solution
time increases by a factor of 4 or 5 when $N$ quadruples, exhibiting
the almost linear complexity of Algorithm \ref{alg:solve}.

The results of the second and third velocity fields are summarized in
Tables \ref{tbl:2DPML2} and \ref{tbl:2DPML3}, respectively. The
quantitative behavior of these tests is similar to the one of the first
velocity field. In all cases, the GMRES iteration converges in about
20 iterations with the sweeping preconditioner.

Secondly, we study how the sweeping preconditioner behaves when $q$
(the number of discretization points per wavelength) varies.  We fix
$\frac{\omega}{2\pi}$ to be $32$ and let $q$ be $8,16,32,64$.  The test
results for the three velocity fields are summarized in Tables
\ref{tbl:2DSPL1}, \ref{tbl:2DSPL2}, and \ref{tbl:2DSPL3}. These
results show that the number of iterations remains to scale at most
logarithmically and the running time of the solution algorithm scales
roughly linearly with respect to the number of unknowns.

\begin{table}[h!]
  \begin{center}
    \begin{tabular}{|cccc|cc|cc|}
      \hline
      \multicolumn{4}{|c|}{} & \multicolumn{2}{|c|}{Test 1} & \multicolumn{2}{|c|}{Test 2}\\
      \hline
      $\omega/(2\pi)$ & $q$ & $N=n^2$ & $T_{\text{setup}}$ & $N_{\text{iter}}$ & $T_{\text{solve}}$ & $N_{\text{iter}}$ & $T_{\text{solve}}$\\
      \hline
      32 & 8  & $255^2$  & 9.19e-01 & 15 & 1.65e+00 & 15 & 1.61e+00\\
      32 & 16 & $511^2$  & 3.91e+00 & 14 & 6.94e+00 & 15 & 7.22e+00\\
      32 & 32 & $1023^2$ & 1.59e+01 & 17 & 8.87e+01 & 17 & 9.39e+01\\
      32 & 64 & $2047^2$ & 6.68e+01 & 19 & 3.74e+02 & 20 & 4.15e+02\\
      \hline
    \end{tabular}
  \end{center}
  \caption{Results of velocity field 1 with varying $q$.}
  \label{tbl:2DSPL1}
\end{table}

\begin{table}[h!]
  \begin{center}
    \begin{tabular}{|cccc|cc|cc|}
      \hline
      \multicolumn{4}{|c|}{} & \multicolumn{2}{|c|}{Test 1} & \multicolumn{2}{|c|}{Test 2}\\
      \hline
      $\omega/(2\pi)$ & $q$ & $N=n^2$ & $T_{\text{setup}}$ & $N_{\text{iter}}$ & $T_{\text{solve}}$ & $N_{\text{iter}}$ & $T_{\text{solve}}$\\
      \hline
      32 & 8  & $255^2$  & 9.28e-01 & 19 & 2.14e+00 & 16 & 1.73e+00\\
      32 & 16 & $511^2$  & 3.69e+00 & 17 & 1.29e+01 & 15 & 1.13e+01\\
      32 & 32 & $1023^2$ & 1.58e+01 & 24 & 1.13e+02 & 15 & 7.16e+01\\
      32 & 64 & $2047^2$ & 6.63e+01 & 26 & 5.29e+02 & 17 & 3.47e+02\\
      \hline
    \end{tabular}
  \end{center}
  \caption{Results of velocity field 2 with varying $q$.}
  \label{tbl:2DSPL2}
\end{table}

\begin{table}[h!]
  \begin{center}
    \begin{tabular}{|cccc|cc|cc|}
      \hline
      \multicolumn{4}{|c|}{} & \multicolumn{2}{|c|}{Test 1} & \multicolumn{2}{|c|}{Test 2}\\
      \hline
      $\omega/(2\pi)$ & $q$ & $N=n^2$ & $T_{\text{setup}}$ & $N_{\text{iter}}$ & $T_{\text{solve}}$ & $N_{\text{iter}}$ & $T_{\text{solve}}$\\
      \hline
      32 & 8  & $255^2$  & 1.00e+00 & 16 & 1.73e+00 & 16 & 1.81e+00\\
      32 & 16 & $511^2$  & 3.66e+00 & 14 & 1.34e+01 & 18 & 1.87e+01\\
      32 & 32 & $1023^2$ & 1.52e+01 & 18 & 8.16e+01 & 19 & 9.22e+01\\
      32 & 64 & $2047^2$ & 6.57e+01 & 19 & 3.99e+02 & 21 & 4.62e+02\\
      \hline
    \end{tabular}
  \end{center}
  \caption{Results of velocity field 3 with varying $q$.}
  \label{tbl:2DSPL3}
\end{table}

Let us compare these numerical results with the ones from the previous
paper \cite{EngquistYing:10a}. The algorithms proposed in this paper
are implemented in Matlab, while the ones in \cite{EngquistYing:10a}
are implemented in C++ with compiler optimization. Hence, direct
comparison of the running time is clearly in favor of the algorithms
in the previous paper. We would expect the running time of the
algorithms in this paper to improve by a factor of 2 to 3 when
implemented in optimized C++ code. Even with this implementational
disadvantage, the setup time $T_{\text{setup}}$ of the current
approach is about twenty times faster. This is mainly due to the fact
that the implementation of the LU factorization is much more efficient
compared to our implementation of the hierarchical matrix algebra in
\cite{EngquistYing:10a}. On the other hand, the number of iterations
$N_{\text{iter}}$ and solution time $T_{\text{solve}}$ of the current
algorithms are higher. This is because in \cite{EngquistYing:10a} $T_m$
are represented by numerical low-rank approximations of the full
half-space Green's function while in the current approach $T_m$ are
approximated based on a modified Helmholtz equation in a truncated
domain.

\subsection{Scattering problem}

The sweeping preconditioner proposed in this paper can also be
extended to scattering problems. Let us consider a simple case where
the scatterer is a sound soft disk centered at the origin with radius
$r_0$. In polar coordinates, the scattered field satisfies the
following equations
\begin{eqnarray*}
  &&\frac{1}{r} \left( r u_r \right)_r + \frac{1}{r^2} u_{\theta\theta} + \frac{w^2}{c^2(r,\theta)} u = f\\
  &&u(r_0,\theta) = -u_{\text{inc}}(r_0,\theta),
\end{eqnarray*}
where $u_{\text{inc}}$ is the incident field and the Sommerfeld
boundary condition is specified for $r$ goes to infinity. One way to
solve this scattering problem is to truncate the domain at $r=r_1$ for
some $r_1>r_0$ and apply the PML condition at $r=r_1$. We can then
apply the sweeping preconditioner in the radial direction from $r=r_1$
to $r=r_0$. In the following example, $c(r,\theta)=1$, $r_0=0.15$ and
$r_1=0.5$. The polar grid is determined so that the each wavelength is
discretized with $q=8$ points. For each fixed $\omega$, two incident
fields are used: one is the Green's function centered at $(-0.2,0.2)$
and the other is the plane wave $\exp(-i \omega x_2)$ traveling
towards the negative $x_2$ direction. We perform tests for
$\frac{w}{2\pi} = 16,32,64,128,256$ and the numerical results are
reported in Table \ref{tbl:2DRAD}.

\begin{table}[h!]
  \begin{center}
    \includegraphics[height=2.3in]{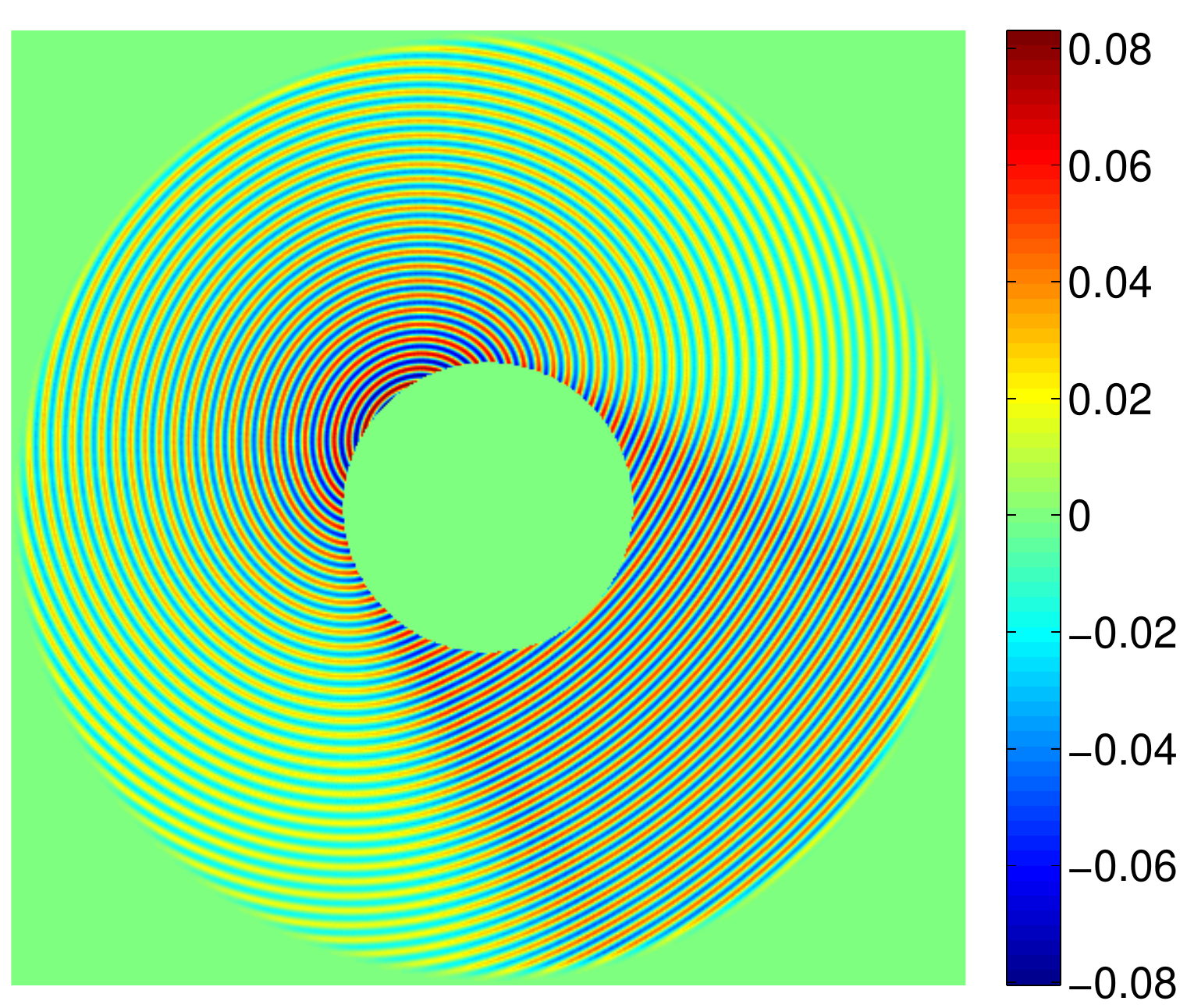}
    \includegraphics[height=2.3in]{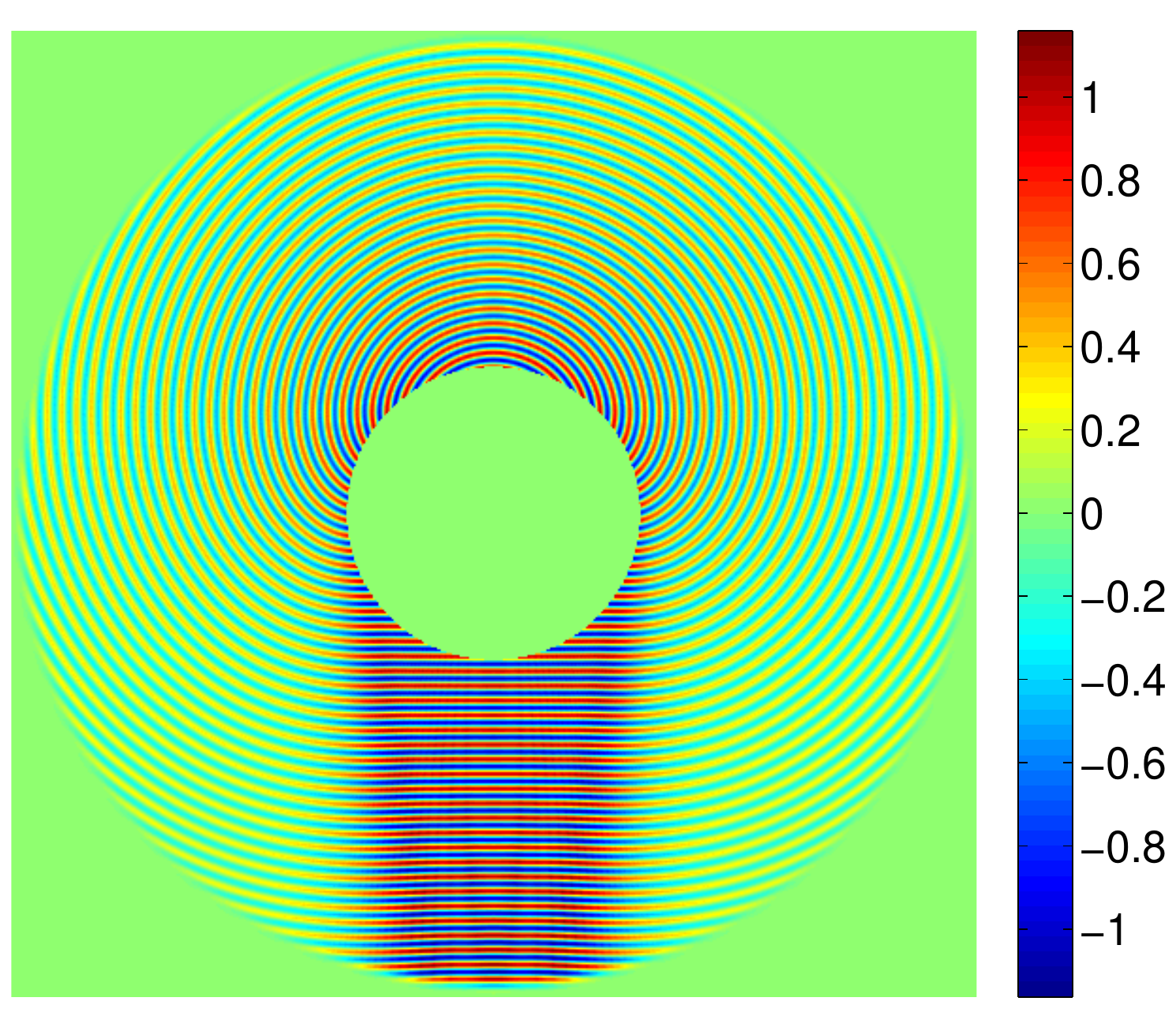}\\
    \begin{tabular}{|cccc|cc|cc|}
      \hline
      \multicolumn{4}{|c|}{} & \multicolumn{2}{|c|}{Incident field 1} & \multicolumn{2}{|c|}{Incident field 2}\\
      \hline
      $\omega/(2\pi)$ & $q$ & $N$ & $T_{\text{setup}}$ & $N_{\text{iter}}$ & $T_{\text{solve}}$ & $N_{\text{iter}}$ & $T_{\text{solve}}$\\
      \hline
      16  & 8 & $45 \times 403$ & 6.11e-01 & 7 & 3.61e-01 & 7 & 2.25e-01\\
      32  & 8 & $90 \times 805$ & 2.61e+00 & 7 & 1.11e+00 & 7 & 1.11e+00\\
      64  & 8 & $180\times 1609$& 1.17e+01 & 7 & 4.92e+00 & 7 & 4.90e+00\\
      128 & 8 & $359\times 3217$& 4.95e+01 & 7 & 2.10e+01 & 7 & 2.06e+01\\
      256 & 8 & $717\times 6434$& 1.99e+02 & 7 & 9.01e+01 & 7 & 9.00e+01\\
      \hline
    \end{tabular}
  \end{center}
  \caption{Results of the scattering problem.
    Top: Scattered fields for two incident waves with $\omega/(2\pi)=64$.
    Bottom: Results for different $\omega$.  }
  \label{tbl:2DRAD}
\end{table}

\section{Preconditioner in 3D}
\label{sec:3Dpre}

The presentation of the 3D preconditioner follows the layout of the 2D
case.

\subsection{Discretization and sweeping factorization}
\label{sec:3Dpredisc}

\begin{figure}[h!]
  \begin{center}
    \includegraphics{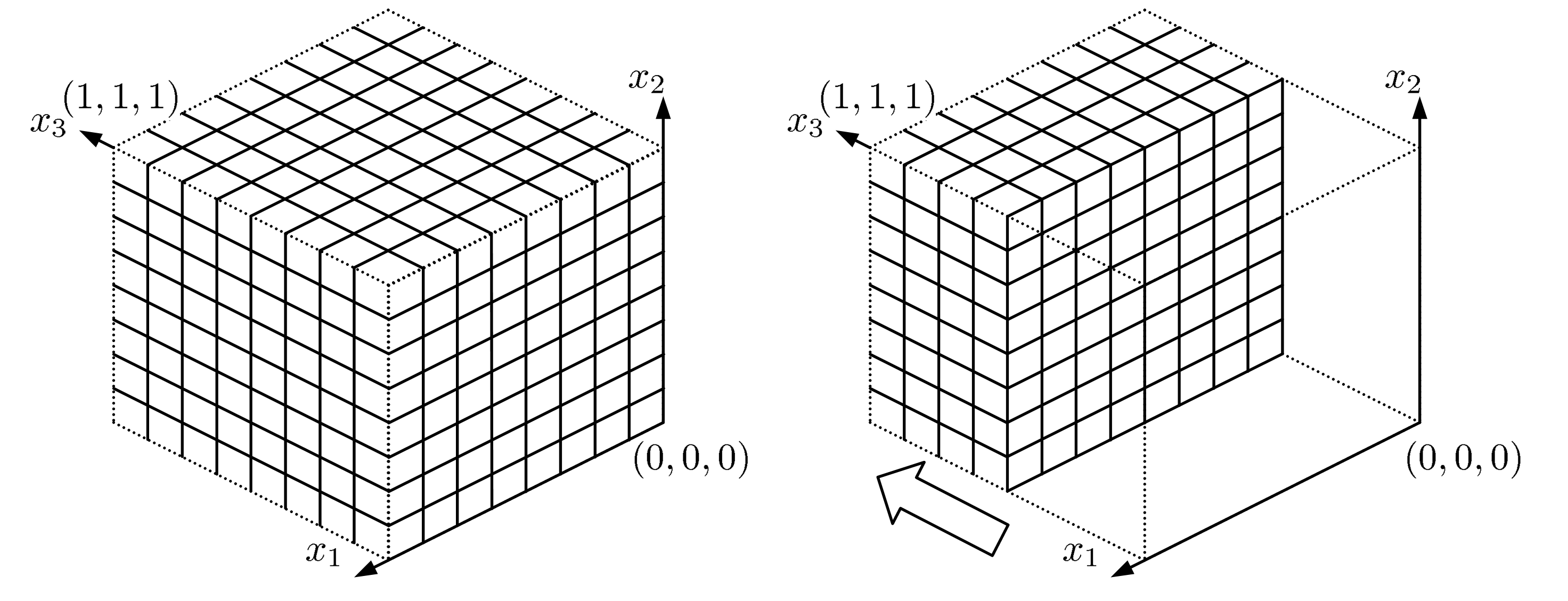}
  \end{center}
  \caption{Left: Discretization grid in 3D. Right: Sweeping order in
    3D. The remaining grid shows the unknowns yet to be processed.}
  \label{fig:3Dgrid}
\end{figure}

The computational domain is $D=(0,1)^3$. Similar to the 2D case,
assume that the Dirichlet boundary condition is used on the side
$x_3=1$ and the Sommerfeld boundary condition is enforced on other
sides. Define
\[
\sigma_1(t) = \sigma_2(t) =
\begin{cases}
  \frac{C}{\eta}\cdot \left( \frac{t-\eta}{\eta} \right)^2 & t \in [0,\eta]\\
  0 & t\in [\eta, 1-\eta] \\
  \frac{C}{\eta}\cdot \left( \frac{t-1+\eta}{\eta} \right)^2 & t \in [1-\eta,1]
\end{cases},
\quad
\sigma_3(t) = 
\begin{cases}
  \frac{C}{\eta}\cdot \left( \frac{t-\eta}{\eta} \right)^2 & t \in [0,\eta]\\
  0 & t\in [\eta, 1]
\end{cases},
\]
and
\[
s_1(x_1) = \left( 1+i \frac{\sigma(x_1)}{\omega} \right)^{-1},\quad
s_2(x_2) = \left( 1+i \frac{\sigma(x_2)}{\omega} \right)^{-1},\quad
s_3(x_3) = \left( 1+i \frac{\sigma(x_3)}{\omega} \right)^{-1}.
\]
The PML approach replaces $\p_1$, $\p_2$, and $\p_2$ with
$s_1(x_1)\p_1$, $s_2(x_2)\p_2$, and $s_3(x_3)\p_3$, respectively. This
effectively provides a damping layer of width $\eta$ near the sides
with Sommerfeld condition. The resulting equation takes the form
\begin{eqnarray*}
  \left( (s_1\p_1)(s_1\p_1) + (s_2\p_2)(s_2\p_2) + (s_3\p_3)(s_3\p_3) + \frac{\omega^2}{c^2(x)} \right) u = f && x\in(0,1)^3,\\
  u = 0 && x \in \p\left([0,1]^3\right).
\end{eqnarray*}
We assume that $f(x)$ is supported inside $[\eta,1-\eta]\times
[\eta,1-\eta]\times [\eta,1]$ (away from the PML). Dividing the
above equation by $s_1(x_1) s_2(x_2) s_3(x_3)$ results
\[
\left( \p_1\left(\frac{s_1}{s_2s_3} \p_1\right)+\p_2\left(\frac{s_2}{s_1s_3} \p_2\right)+\p_3\left(\frac{s_3}{s_1s_2} \p_3\right) + 
  \frac{\omega^2}{s_1s_2s_3c^2(x)} \right) u = f.
\]
The domain $[0,1]^3$ is discretized with a Cartesian grid with spacing
$h = 1/(n+1)$, where the number of points $n$ in each dimension is
proportional to $\omega$. The interior points of this grid are
\[
\P = \{ p_{i,j,k} = \left(ih,jh,kh\right) : 1\le i,j,k \le  n\}
\]
(see Figure \ref{fig:3Dgrid} (left)) and the total number of grid
points is $N=n^3$.

We denote by $u_{i,j,k}$, $f_{i,j,k}$, and $c_{i,j,k}$ the values of
$u(x)$, $f(x)$, and $c(x)$ at point $x_{i,j,k} = (ih,jh,kh)$. The
standard $7$-point stencil finite difference method writes down the
equation at points in $\P$ using central difference. The resulting
equation at $(ih,jh,kh)$ is
\begin{multline}
  \frac{1}{h^2} \left(\frac{s_1}{s_2s_3}\right)_{i-\hf,j,k} u_{i-1,j,k}
  + \frac{1}{h^2} \left(\frac{s_1}{s_2s_3}\right)_{i+\hf,j,k} u_{i+1,j,k}
  + \frac{1}{h^2} \left(\frac{s_2}{s_1s_3}\right)_{i,j-\hf,k} u_{i,j-1,k} \\
  + \frac{1}{h^2} \left(\frac{s_2}{s_1s_3}\right)_{i,j+\hf,k} u_{i,j+1,k}
  + \frac{1}{h^2} \left(\frac{s_3}{s_1s_2}\right)_{i,j,k-\hf} u_{i,j,k-1}
  + \frac{1}{h^2} \left(\frac{s_3}{s_1s_2}\right)_{i,j,j+\hf} u_{i,j,k+1}\\
  + \left(
    \frac{\omega^2}{(s_1s_2s_3)_{i,j,k}\cdot c_{i,j,k}^2 }
    - \left(\cdots
    \right) 
  \right)
  u_{i,j,k} = f_{i,j,k}
  \label{eq:helmd3}
\end{multline} 
with $u_{i',j',k'}$ equal to zero for $(i',j',k')$ that violates $1\le
i',j',k'\le n$. Here $(\cdots)$ stands for the sum of the six
coefficients appeared in the first two lines. We order $u_{i,j,k}$ and
$f_{i,j,k}$ by going through the three dimensions in order and denote
the vectors containing by
\begin{eqnarray*}
&& u = \left(u_{1,1,1}, u_{2,1,1},\ldots,u_{n,1,1}, \ldots,u_{1,n,n}, u_{2,n,n},\ldots,u_{n,n,n}\right)^t\\
&& f = \left(f_{1,1,1}, f_{2,1,1},\ldots,f_{n,1,1}, \ldots,f_{1,n,n}, f_{2,n,n},\ldots,f_{n,n,n}\right)^t
\end{eqnarray*}
The discrete system of \eqref{eq:helmd3} takes the form $A u = f$. We
further introduce a block version.  Define $\P_m$ to be the indices in
the $m$-th row
\[
\P_m = \{p_{1,1,m}, p_{2,1,m}, \ldots, p_{n,n,m} \}
\]
and introduce
\[
u_m = \left( u_{1,1,m}, u_{2,1,m},\ldots,u_{n,n,m} \right)^t\quad\text{and}\quad
f_m = \left( f_{1,1,m}, f_{2,1,m},\ldots,f_{n,n,m} \right)^t.
\]
Then
\[
u = (u_1^t, u_2^t, \ldots, u_n^t)^t,\quad
f = (f_1^t, f_2^t, \ldots, f_n^t)^t.
\]
Using these notations, the system $Au=f$ takes the following block
tridiagonal form
\[
\begin{pmatrix}
  A_{1,1} & A_{1,2} & & \\
  A_{2,1} & A_{2,2} & \ddots & \\
  & \ddots & \ddots & A_{n-1,n}\\
  & & A_{n,n-1} & A_{n,n}
\end{pmatrix}
\begin{pmatrix}
  u_1\\
  u_2\\
  \vdots\\
  u_n
\end{pmatrix}
=
\begin{pmatrix}
  f_1\\
  f_2\\
  \vdots\\
  f_n
\end{pmatrix}
\]
where each block is of size $n^2\times n^2$ and the off-diagonal
blocks are diagonal matrices. The sweeping factorization takes the
same form as the 2D one \eqref{eq:Afact}. In order to design an
efficient preconditioner, the main task task is to construct
approximations for the Schur complement matrix $T_m: g_m \rightarrow
v_m$, which maps an external force $g_m$ loaded only on the $m$-th
layer to the solution $v_m$ restricted to the same layer. Following
the central idea of pushing the PML right next to $x_3=mh$, we define
\[
s^m_3(x_3) = \left( 1+i\frac{\sigma_3(x_3-(m-b)h)}{\omega} \right)^{-1}.
\]
and introduce an auxiliary problem on the domain $D_m = [0,1]\times
[0,1]\times [(m-b)h,(m+1)h]$:
\begin{eqnarray}
\left( (s_1\p_1)(s_1\p_1) + (s_2\p_2)(s_2\p_2) + (s^m_3\p_3)(s^m_3\p_3) + \frac{\omega^2}{c^2(x)} \right) u = f && x\in D_m, \label{eq:helma3}\\
u = 0 && x \in \p D_m . \nonumber
\end{eqnarray}
This equation is then discretized with the subgrid
\[
\G_m = \{p_{i,j,k}, 1\le i,j \le n, m-b+1 \le k \le m \}
\]
of the original grid $\P$. The resulting $bn^2 \times bn^2$ discrete
Helmholtz operator is denoted by $H_m$. The operator $\wt{T}_m: g_m
\rightarrow v_m$ defined by
\[
\begin{pmatrix}
  *\\ \vdots \\ * \\ v_m
\end{pmatrix}
\approx 
H_m^{-1}
\begin{pmatrix}
  0\\ \vdots \\ 0 \\ g_m
\end{pmatrix}
\]
is an approximation of the Schur complement matrix $T_m$.  Since $H_m$
comes from the 7-point stencil with $b$ layers, this can be viewed as
a quasi-2D problem, which can be solved efficiently using a modified
version of the multifrontal method
\cite{DuffReid:83,George:73,Liu:92}.

\begin{figure}[h!]
  \begin{center}
    \includegraphics{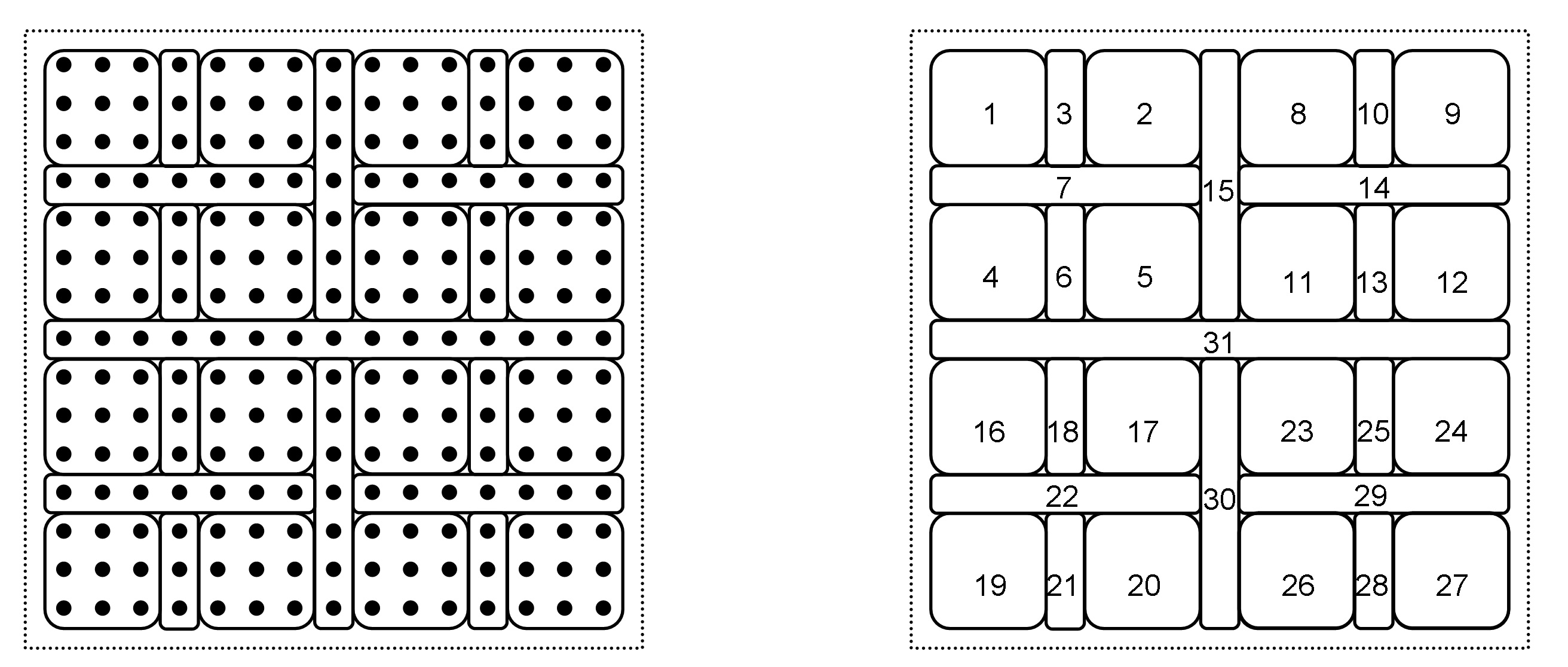}
  \end{center}
  \caption{Multifrontal algorithm on a $15 \times 15$ two dimensional
    Cartesian grid. Left: The unknowns are clustered into groups
    hierarchically to minimizes the boundary between different
    groups. Right: Elimination order of different groups. The groups
    are eliminated in the increasing order of their indices.}
  \label{fig:multifrontal}
\end{figure}

The main idea of the multifrontal method is simple yet elegant. Take a
$n\times n$ 2D grid as an example and use $M$ to denote the discrete
operator resulted from a local stencil. The multifrontal method
reorders the unknowns hierarchically in order to minimize the fill-ins
of the $LDL^t$ factorization of $M$. For the $n \times n$ Cartesian
grid, one possible ordering is given in Figure \ref{fig:multifrontal}
where the unknowns are clustered into groups and the groups are
ordered hierarchically. The construction of the $LDL^t$ factorization
eliminates the unknowns group by group. The dominating cost of the
algorithm is spent in inverting the unknowns of the last few groups
and the overall cost is $O(n^3)$, cubic in terms of the size of the
last group. Moreover, the $L$ matrix is never constructed explicitly
in the multifrontal method. Instead it is stored and applied as a
sequence of (block) row operations for the sake of
efficiency. Applying $M^{-1}$ to an arbitrary vector using the result
of the multifrontal algorithm takes $O(n^2\log n)$ steps. In the
current setting, we adopt the same hierarchical partitioning in the
$(x_1,x_2)$ plane, while keeping the unknowns with the same $x_1$ and
$x_2$ indices in the same group. Since now the size of the last group
is of order $O(bn)$, the construction phase of the multifrontal method
takes $O(b^3 n^3)$ steps and applying to an arbitrary vector takes
$O(b^2 n^2 \log n)$ steps.

\subsection{Approximate inversion and preconditioner}

Let us now combine the multifrontal method into Algorithms
\ref{alg:setupext} and \ref{alg:solveext} to build the approximate
inverse of $H$. Similar to the 2D case, we define 
\[
u_F = (u_1^t,\ldots,u_{b}^t)^t \quad
f_F = (f_1^t,\ldots,f_{b}^t)^t.
\]
and write
\[
\begin{pmatrix}
  A_{F,F} & A_{F,b+1} & & \\
  A_{b+1,F} & A_{b+1,b+1} & \ddots & \\
  & \ddots & \ddots & A_{n-1,n}\\
  & & A_{n,n-1} & A_{n,n}
\end{pmatrix}
\begin{pmatrix}
  u_F\\
  u_{b+1}\\
  \vdots\\
  u_n
\end{pmatrix}
=
\begin{pmatrix}
  f_F\\
  f_{b+1}\\
  \vdots\\
  f_n
\end{pmatrix}.
\]
The goal of the construction of the approximate sweeping factorization
of $A$ is to compute $\wt{T}_m$ and the algorithm consists of the
following steps.
\begin{algo}
  Construction of the approximate sweeping factorization of $A$ with
  moving PML.
  \label{alg:setup3}
\end{algo}
\begin{algorithmic}[1]
  \STATE Let $\G_F$ be the subgrid of the first $b$ layers and
  $H_F = A_{F,F}$. Construct the multifrontal factorization of $H_F$
  by partitioning $\G_F$ hierarchically in the $(x_1,x_2)$ plane.
  
  \FOR{$m=b+1,\ldots,n$}

  \STATE Let $\G_m = \{p_{i,j,k}, 1 \le i,j\le n, m-b+1 \le k \le m\}$
  and $H_m$ be the system of \eqref{eq:helma3} on $\G_m$.  Construct
  the multifrontal factorization of $H_m$ by partitioning $\G_m$
  hierarchically in the $(x_1,x_2)$ plane.
  
  \ENDFOR
\end{algorithmic}
The cost of Algorithm \ref{alg:setup3} is $O(b^3 n^4) = O(b^3
N^{4/3})$. The computation of $u$ from this sweeping factorization is
summarized in the following algorithm
\begin{algo}
  Computation of $u \approx A^{-1}f$ using the sweeping factorization
  of $A$ with moving PML.
  \label{alg:solve3}
\end{algo}
\begin{algorithmic}[1]
  \STATE $u_F = f_F$ and $u_m = f_m$ for $m=b+1,\ldots,n$.
  
  \STATE $u_{b+1} = u_{b+1} - A_{b+1,F} (\wt{T}_F u_F)$. $\wt{T}_F
  u_F$ is computed using the multifrontal factorization of $H_F$.
  
  \FOR{$m=b+1,\ldots,n-1$}
  
  \STATE $u_{m+1} = u_{m+1} - A_{m+1,m} (\wt{T}_m u_m)$. The
  application of $\wt{T}_m u_m$ is done by forming the vector
  $(0,\ldots,0,u_m^t)^t$, applying $H_m^{-1}$ to it using the
  multifrontal factorization of $H_m$, and extracting the value on the
  last layer.
  
  \ENDFOR
  
  \STATE $u_F = \wt{T}_F u_F$. See the previous steps for the
  application of $\wt{T}_F$.

  \FOR{$m=b+1,\ldots,n$}
  
  \STATE $u_m = \wt{T}_m u_m$. See the previous steps for the
  application of $\wt{T}_m$.
  
  \ENDFOR
  
  \FOR{$m=n-1,\ldots,b+1$}
  
  \STATE $u_m = u_m - \wt{T}_m (A_{m,m+1} u_{m+1})$. See the previous
  steps for the application of $\wt{T}_m$.
  
  \ENDFOR

  \STATE $u_F = u_F - \wt{T}_F (A_{F,b+1} u_{b+1})$. See the previous
  steps for the application of $\wt{T}_F$.
  
\end{algorithmic}
The cost of Algorithm \ref{alg:solve3} is $O(b^2 n^3 \log n) = O(b^2 N
\log N)$.

For the stability reason mentioned in Section \ref{sec:2Dpre}, we
apply Algorithms \ref{alg:setup3} and \ref{alg:solve3} to the discrete
operator $A_\alpha$ of the modified system
\[
\Lapl u(x) + \frac{(\omega+i\alpha)^2}{c^2(x)} u(x) = f(x),
\]
where $\alpha$ is an $O(1)$ positive constant. We denote by $M_\alpha:
f \rightarrow u$ the operator defined by Algorithm \ref{alg:solve} for
this modified equation. Since $A_\alpha$ is close to $A$ when $\alpha$
is small, we propose to solve the preconditioner system
\[
M_\alpha A u = M_\alpha f
\]
using the GMRES solver \cite{Saad:03,SaadSchultz:86}. Because the cost
of applying $M_\alpha$ to any vector is $O(N\log N)$, the total cost
of the GMRES solver is $O(N_I N\log N)$, where $N_I$ is the number of
iterations required. As the numerical results in Section
\ref{sec:3Dnum} demonstrate, $N_I$ is essentially independent of the
number of unknowns $N$, thus resulting an algorithm with almost linear
complexity.

\section{Numerical Results in 3D}
\label{sec:3Dnum}

In this section, we present several numerical results to illustrate
the properties of the algorithm described in Section \ref{sec:3Dpre}.
We use GMRES as the iterative solver with relative residue tolerance
equal to $10^{-3}$. 

\begin{figure}[h!]
  \begin{center}
    \begin{tabular}{ccc}
      \includegraphics[height=1.6in]{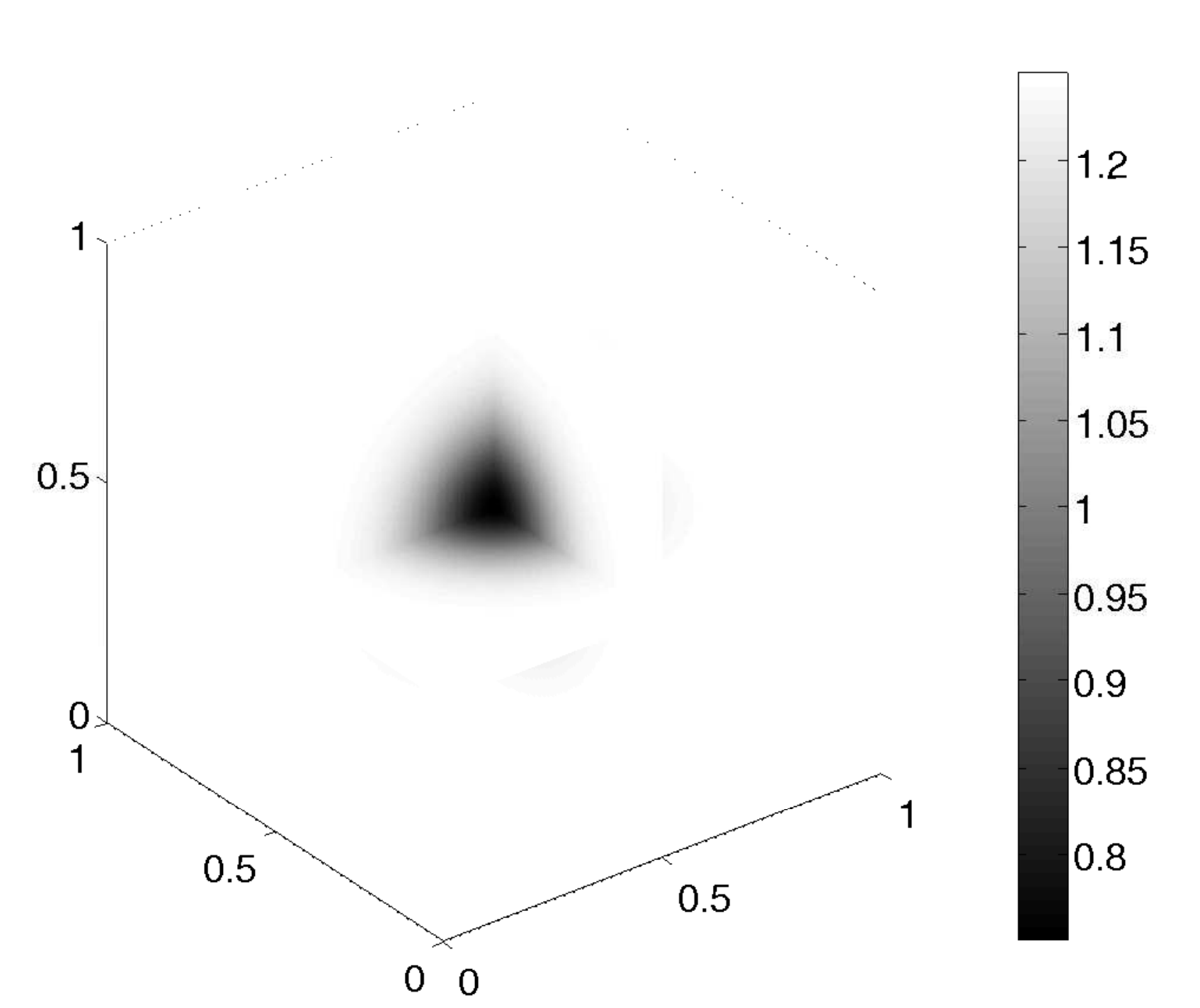}&\includegraphics[height=1.6in]{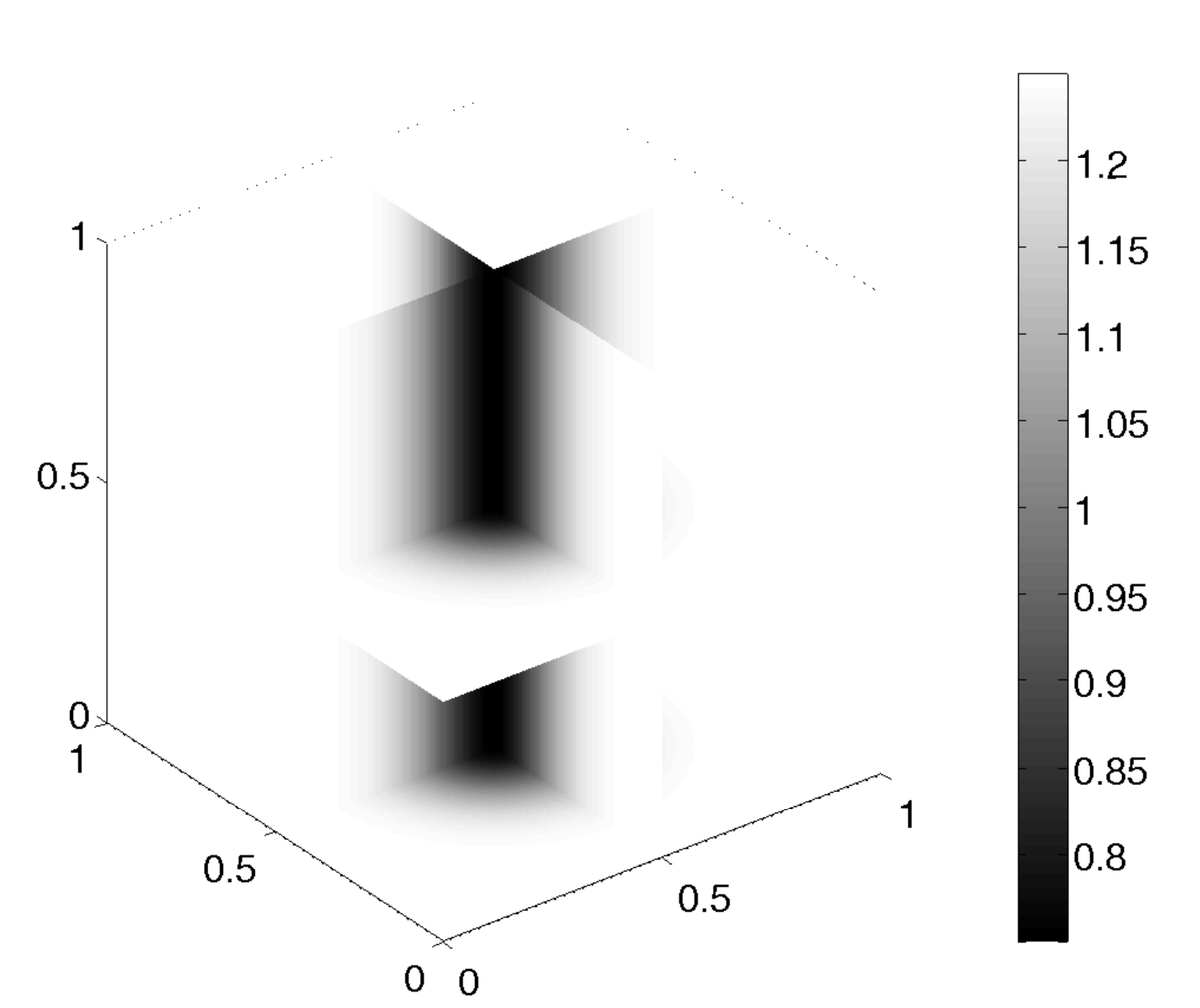}& \includegraphics[height=1.6in]{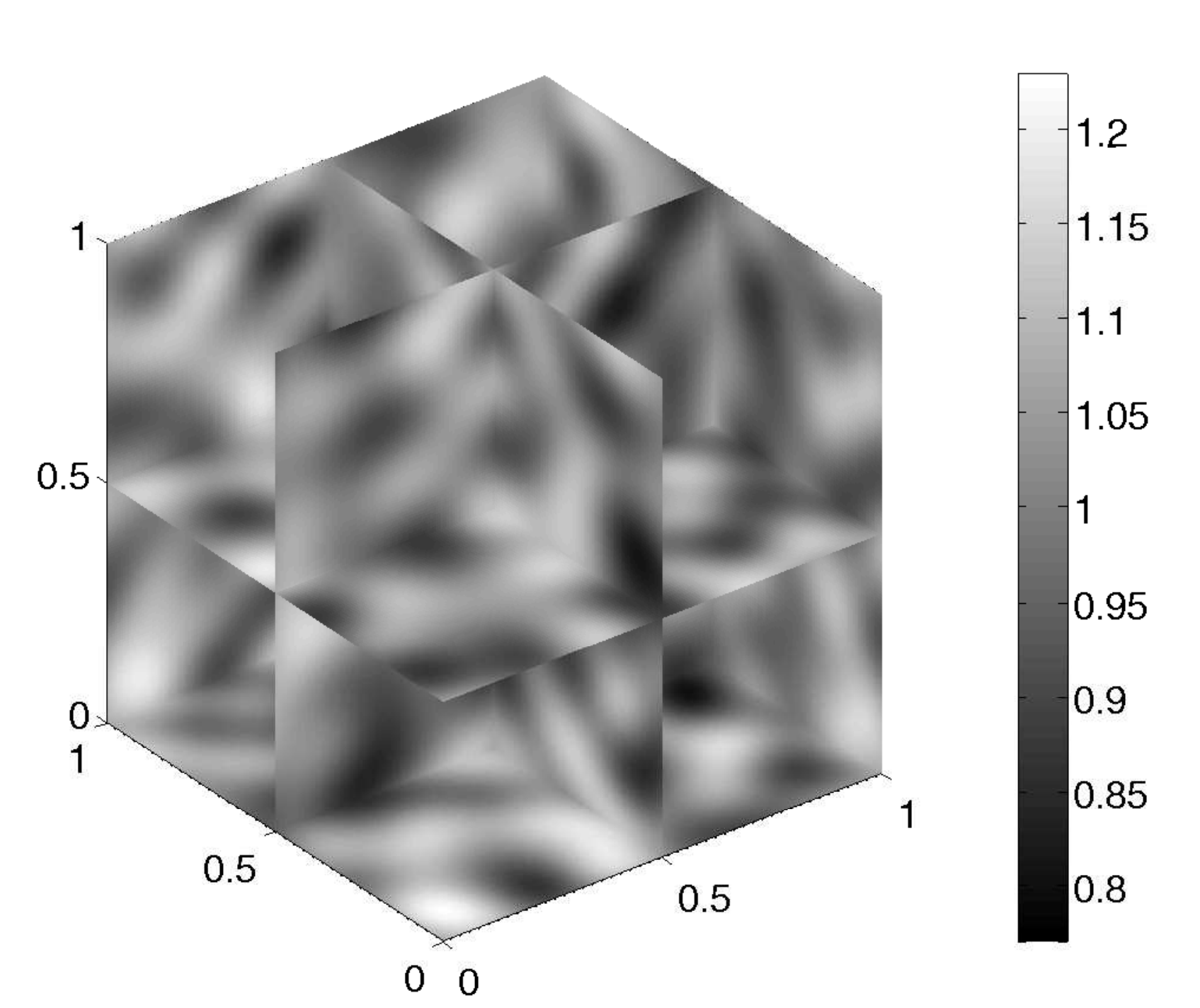}\\
      (a) & (b) & (c)
    \end{tabular}
  \end{center}
  \caption{Test velocity fields.}
  \label{fig:3Dnumspeed}
\end{figure}

The examples in this seciton have the PML boundary condition specified
at all sides. We consider three velocity fields in the domain
$D=(0,1)^3$:
\begin{enumerate}
\item The first velocity field is a converging lens with a Gaussian
  profile at the center of the domain (see Figure
  \ref{fig:3Dnumspeed}(a)).
\item The second velocity field is a vertical waveguide with Gaussian
  cross section (see Figure \ref{fig:3Dnumspeed}(b)).
\item The third velocity field is a random velocity field  (see
  Figure \ref{fig:3Dnumspeed}(c)).
\end{enumerate}

For each problem, we test with two external forces $f(x)$.
\begin{enumerate}
\item The first external force $f(x)$ is a Gaussian point source
  located at $(x_1,x_2,x_3) = (0.5, 0.5, 0.25)$. The response of this
  forcing term generates circular waves propagating at all
  directions. Due to the variations of the velocity field, the
  circular waves are going to bend and form caustics.
\item The second external force $f(x)$ is a Gaussian wave packet whose
  wavelength is comparable to the typical wavelength of the
  domain. This packet centers at $(x_1,x_2,x_3) = (0.5, 0.25, 0.25)$
  and points to the $(0,1,1)$ direction. The response of this forcing
  term generates a Gaussian beam initially pointing towards the
  $(0,1,1)$ direction. 
\end{enumerate}

Firstly, we study how the sweeping preconditioner behaves when
$\omega$ varies. For each velocity field, we perform tests for
$\frac{\omega}{2\pi}$ equal to $5,10,20$. In these tests, we
discretize each wavelength with $q=8$ points and the number of samples
in each dimension is $n=39,79,159$. The $\alpha$ value of the modified
system is set to be equal to 1. The width of the PML is equal to $6h$
(i.e., $b=6$) and the number of layers processed within each iteration
of Algorithms \ref{alg:setup3} and \ref{alg:solve3} is equal to 3
(i.e., $d=3$). The preconditioner sweeps the domain with two fronts
that start from $x_3=0$ and $x_3=1$.

The results of the first velocity field is reported in Table
\ref{tbl:3DPML1}. The two plots show the solutions of the two right
sides on a plane near $x_1 = 0.5$. $T_{\text{setup}}$ is the time used
to construct the preconditioner in seconds.  $N_{\text{iter}}$ is the
number of iterations of the preconditioned GMRES iteration and
$T_{\text{solve}}$ is the solution time. The estimate in Section
\ref{sec:3Dpre} section shows that the setup time scales like
$O(N^{4/3})$. So when $\omega$ doubles, $N$ increases by a factor of
$4$ and $T_{\text{setup}}$ should increase by a factor of $16$. The
numerical results show that the actual growth factor is even lower. A
remarkable feature of the sweeping preconditioner is that in all cases
the preconditioned GMRES solver converges in at most 12
iterations. Finally, we would like to point out that our algorithm is
quite efficient: for the case with $\omega/(2\pi)=20$ with more than
four million unknowns, the solution time is less than 600 seconds.
The results of the second and the third velocity fields are reported
in Tables \ref{tbl:3DPML2} and \ref{tbl:3DPML3}, respectively. In all
tests, the GMRES iteration converges at most 13 iterations when
combined with the new sweeping preconditioner.

\begin{table}[h!]
  \begin{center}
    \includegraphics[height=2.3in]{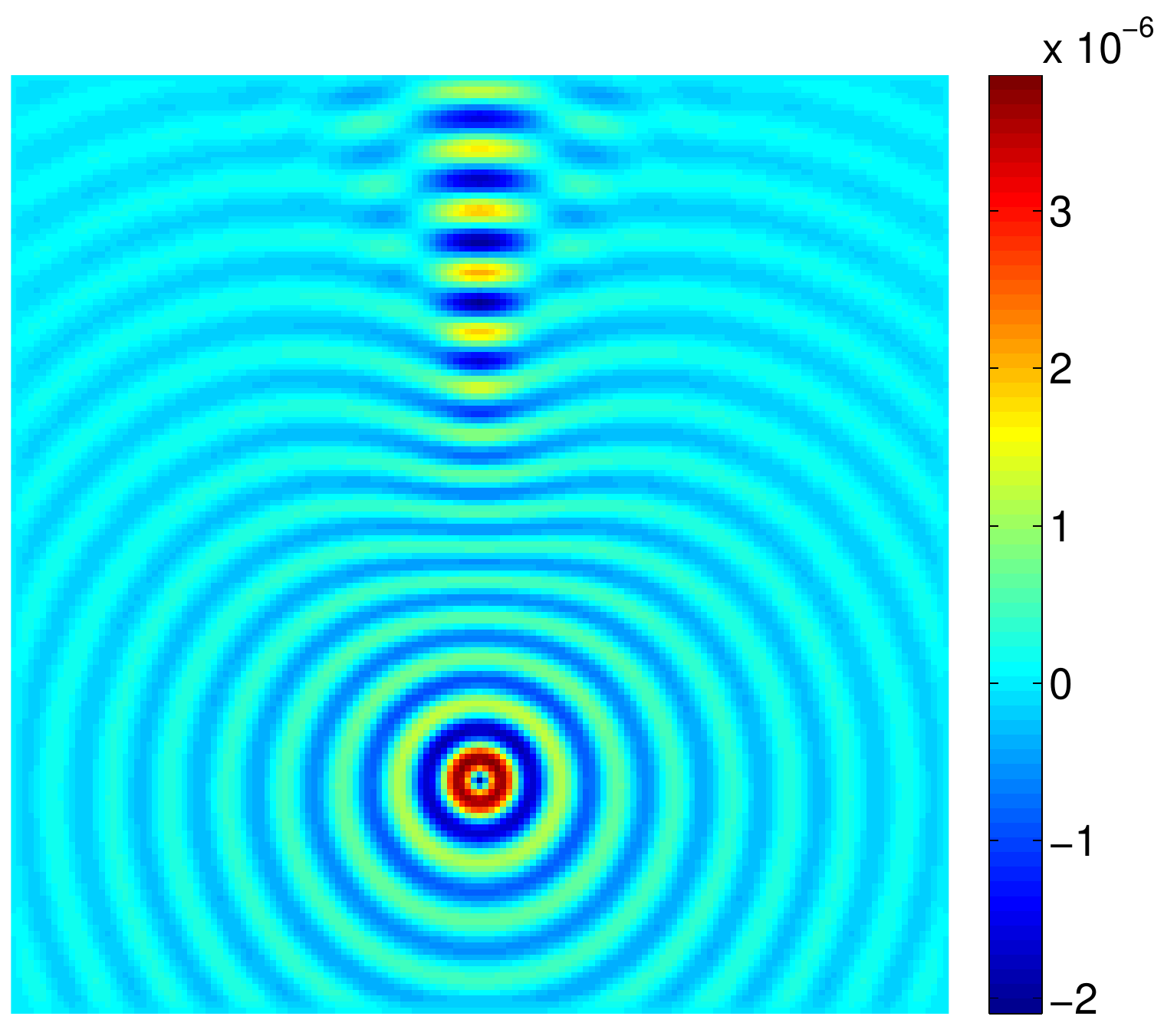}
    \includegraphics[height=2.3in]{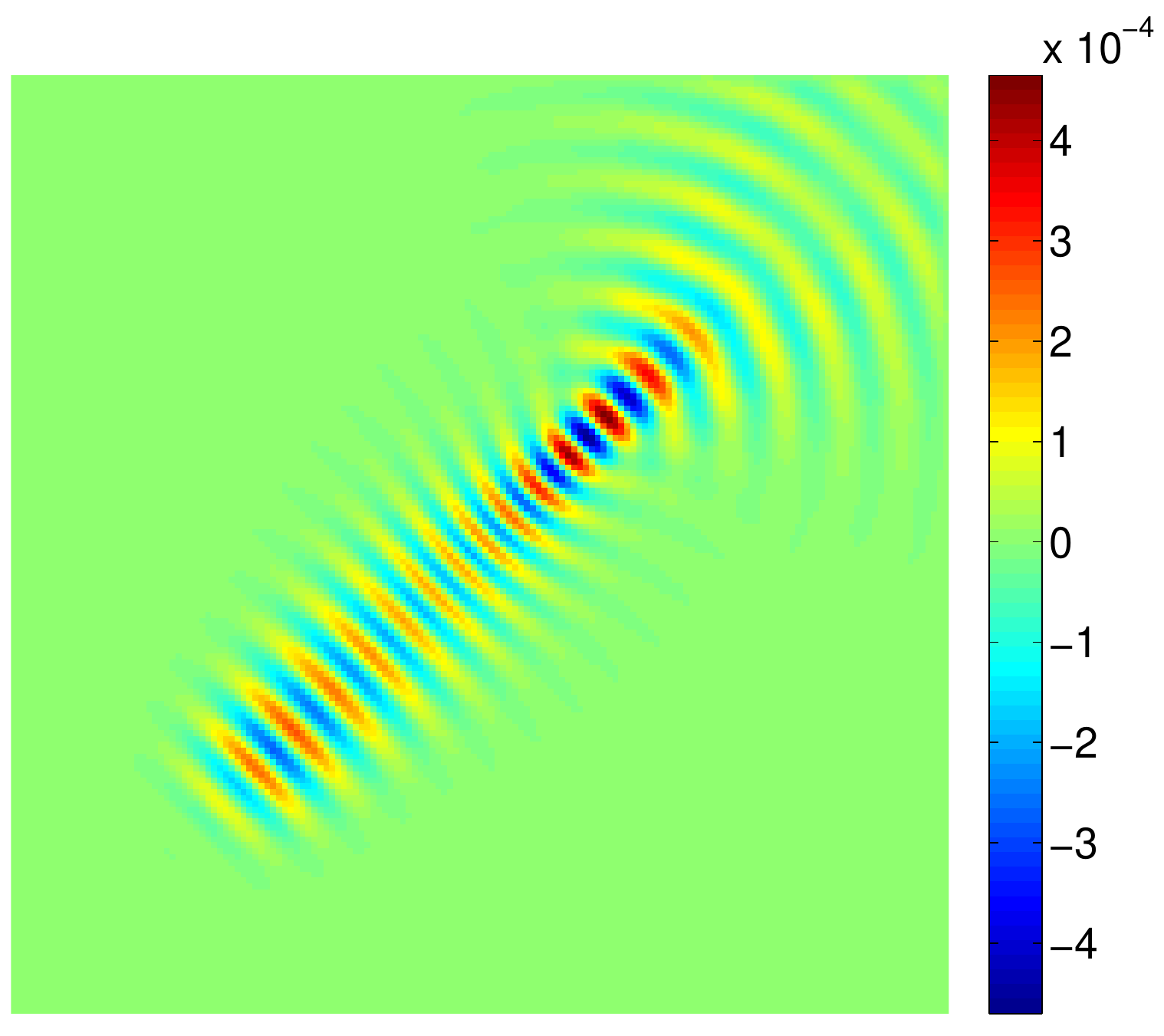}\\
    \begin{tabular}{|cccc|cc|cc|}
      \hline
      \multicolumn{4}{|c|}{} & \multicolumn{2}{|c|}{Test 1} & \multicolumn{2}{|c|}{Test 2}\\
      \hline
      $\omega/(2\pi)$ & $q$ & $N=n^3$ & $T_{\text{setup}}$ & $N_{\text{iter}}$ & $T_{\text{solve}}$ & $N_{\text{iter}}$ & $T_{\text{solve}}$\\
      \hline
      5  & 8 & $39^3$  & 4.80e+00 & 11 & 4.53e+00 & 11 & 4.63e+00\\
      10 & 8 & $79^3$  & 6.37e+01 & 11 & 4.92e+01 & 11 & 4.93e+01\\
      20 & 8 & $159^3$ & 8.27e+02 & 12 & 5.53e+02 & 12 & 5.94e+02\\
      \hline
    \end{tabular}
  \end{center}
  \caption{Results of velocity field 1 with varying $\omega$. 
    Top: Solutions for two external forces with $\omega/(2\pi)=16$ on a plane near $x_1=0.5$.
    Bottom: Results for different $\omega$.  }
  \label{tbl:3DPML1}
\end{table}

\begin{table}[h!]
  \begin{center}
    \includegraphics[height=2.3in]{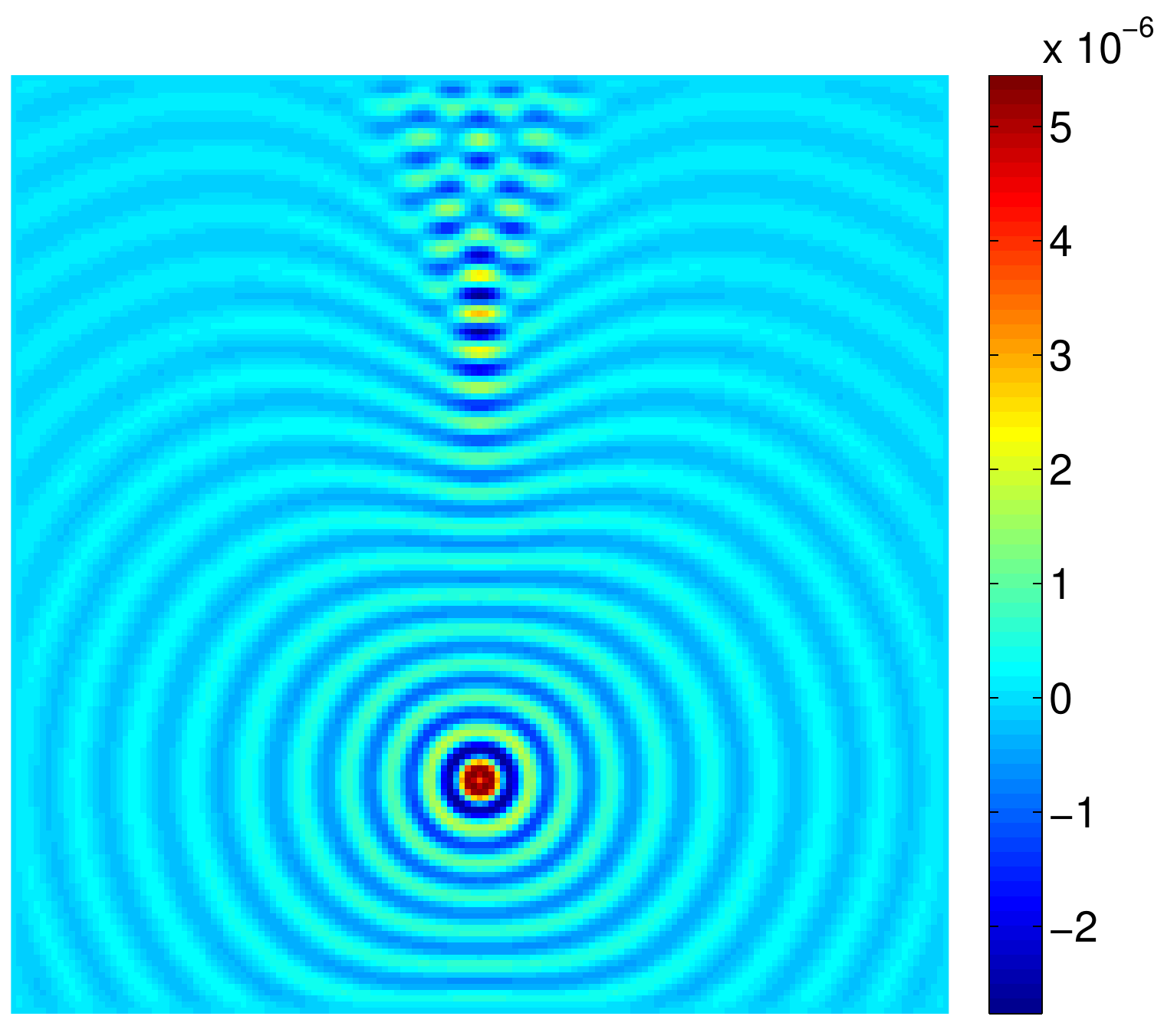}
    \includegraphics[height=2.3in]{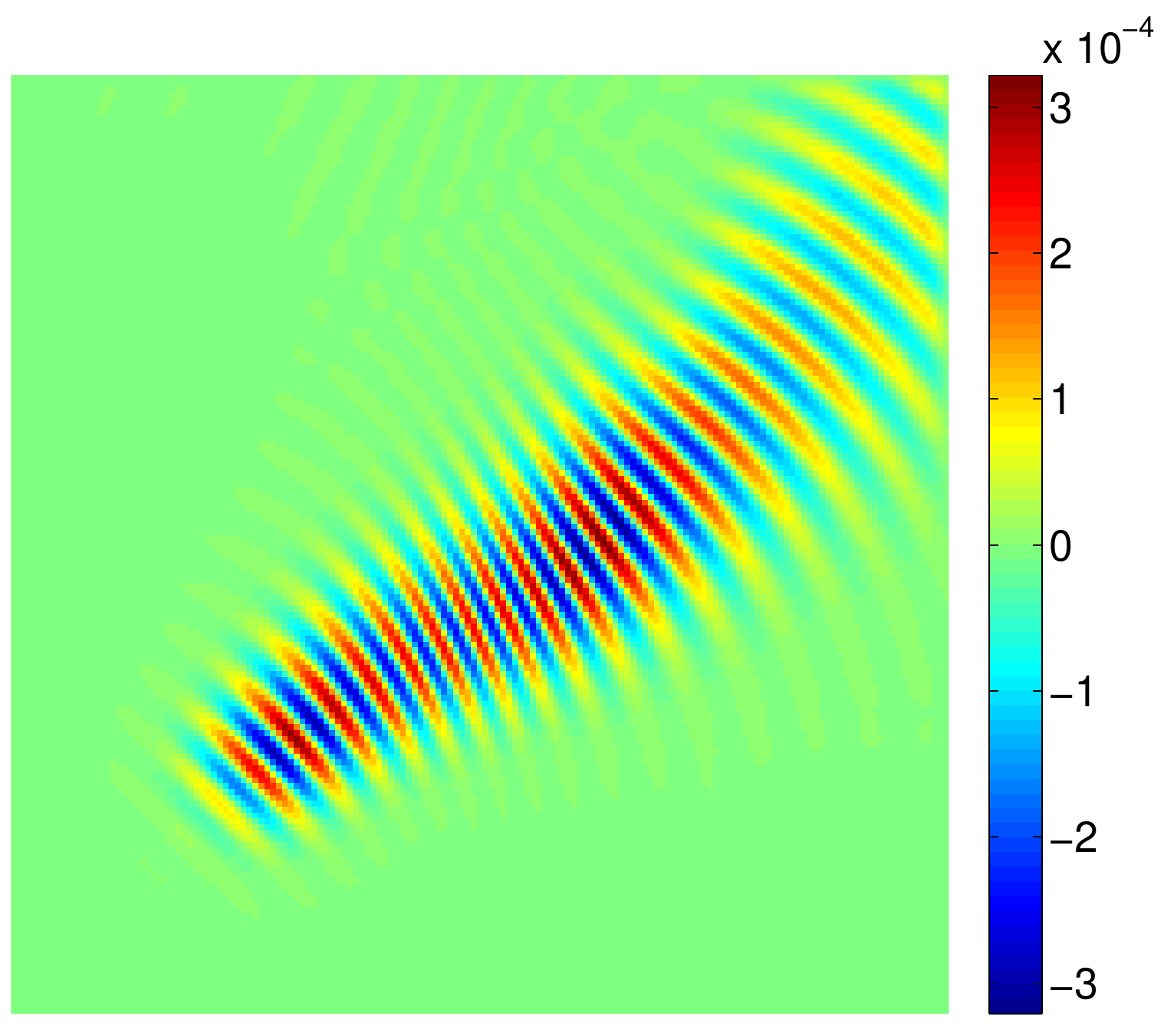}\\
    \begin{tabular}{|cccc|cc|cc|}
      \hline
      \multicolumn{4}{|c|}{} & \multicolumn{2}{|c|}{Test 1} & \multicolumn{2}{|c|}{Test 2}\\
      \hline
      $\omega/(2\pi)$ & $q$ & $N=n^3$ & $T_{\text{setup}}$ & $N_{\text{iter}}$ & $T_{\text{solve}}$ & $N_{\text{iter}}$ & $T_{\text{solve}}$\\
      \hline
      5  & 8 & $39^3$  & 4.83e+00 & 12 & 5.14e+00 & 12 & 5.03e+00\\
      10 & 8 & $79^3$  & 6.76e+01 & 13 & 5.70e+01 & 12 & 5.64e+01\\
      20 & 8 & $159^3$ & 8.24e+02 & 14 & 6.32e+02 & 11 & 5.40e+02\\
      \hline
    \end{tabular}
  \end{center}
  \caption{Results of velocity field 2 with varying $\omega$. 
    Top: Solutions for two external forces with $\omega/(2\pi)=16$ on a plane near $x_1=0.5$.
    Bottom: Results for different $\omega$.  }
  \label{tbl:3DPML2}
\end{table}

\begin{table}[h!]
  \begin{center}
    \includegraphics[height=2.3in]{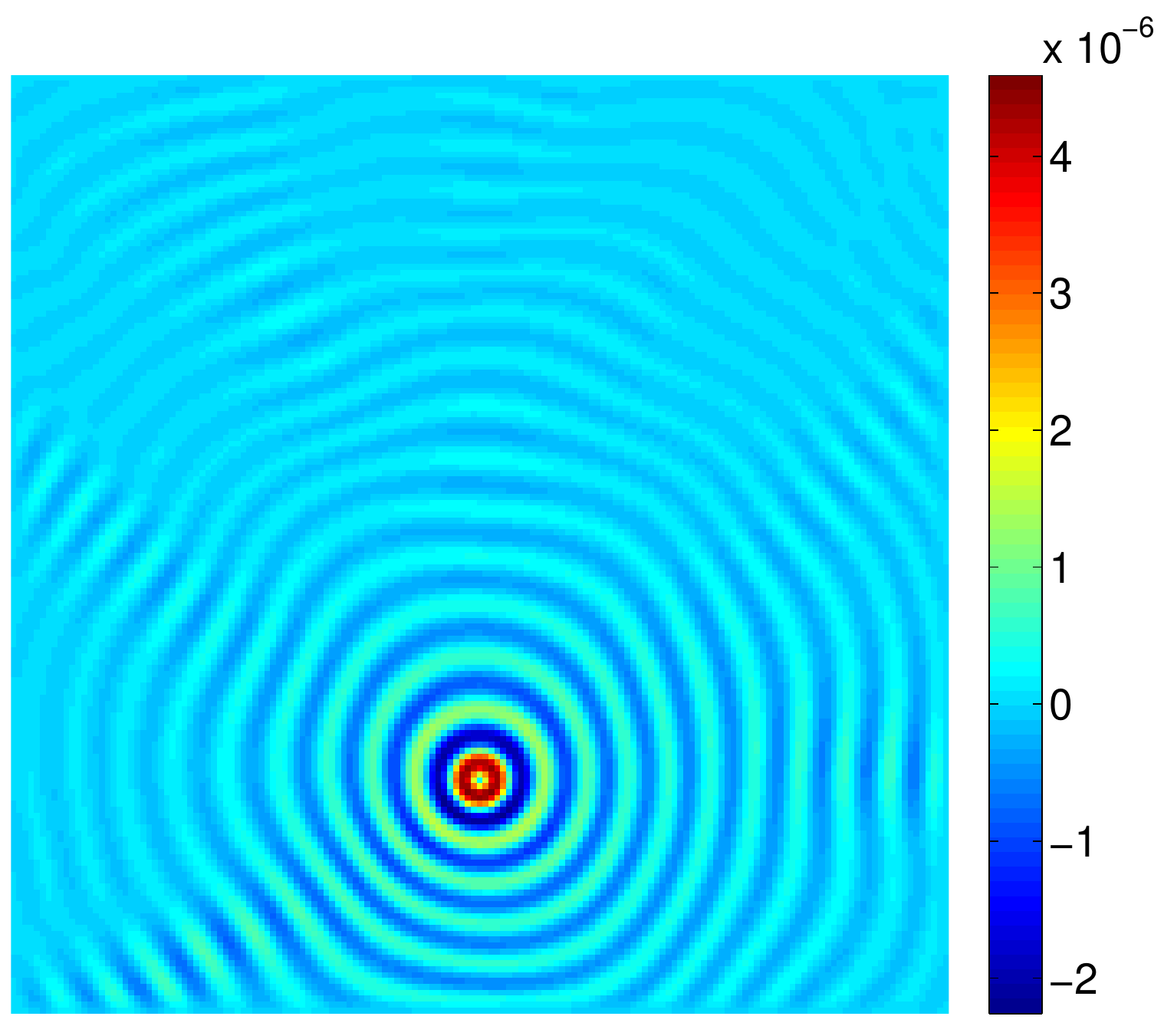}
    \includegraphics[height=2.3in]{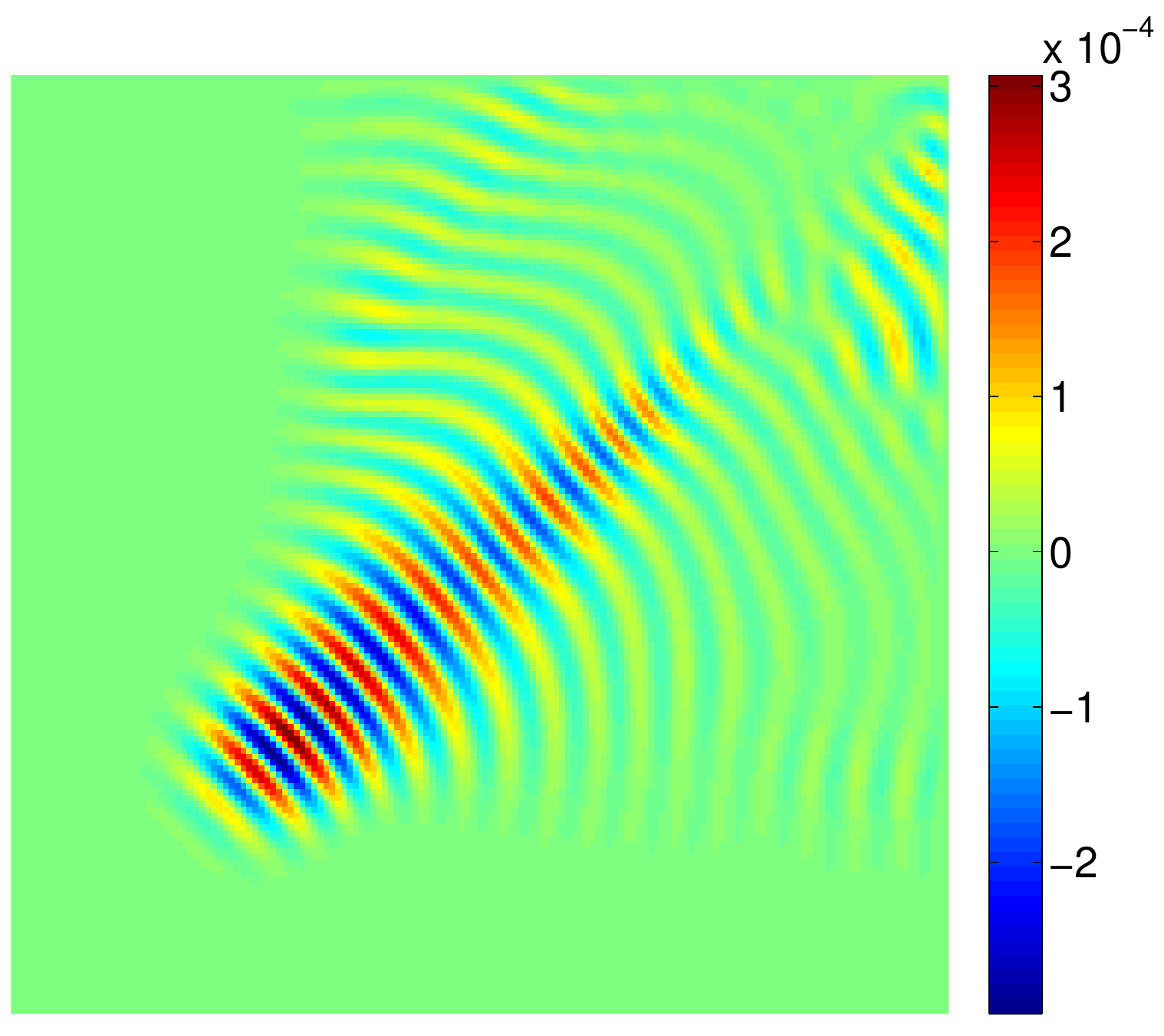}\\
    \begin{tabular}{|cccc|cc|cc|}
      \hline
      \multicolumn{4}{|c|}{} & \multicolumn{2}{|c|}{Test 1} & \multicolumn{2}{|c|}{Test 2}\\
      \hline
      $\omega/(2\pi)$ & $q$ & $N=n^3$ & $T_{\text{setup}}$ & $N_{\text{iter}}$ & $T_{\text{solve}}$ & $N_{\text{iter}}$ & $T_{\text{solve}}$\\
      \hline
      5  & 8 & $39^3$  & 4.85e+00 & 12 & 5.26e+00 & 12 & 5.44e+00\\
      10 & 8 & $79^3$  & 6.69e+01 & 11 & 5.10e+01 & 13 & 5.99e+01\\
      20 & 8 & $159^3$ & 8.42e+02 & 11 & 5.58e+02 & 13 & 6.28e+02\\
      \hline
    \end{tabular}
  \end{center}
  \caption{Results of velocity field 3 with varying $\omega$. 
    Top: Solutions for two external forces with $\omega/(2\pi)=16$ on a plane near $x_1=0.5$.
    Bottom: Results for different $\omega$.  }
  \label{tbl:3DPML3}
\end{table}

Secondly, we study how the sweeping preconditioner behaves when $q$
(the number of discretization points per wavelength) varies.  We fix
$\frac{\omega}{2\pi}$ to be $5$ and let $q$ be $8,16,32$.  The test
results for the three velocity fields are summarized in Tables
\ref{tbl:3DSPL1}, \ref{tbl:3DSPL2}, and \ref{tbl:3DSPL3},
respectively. These results show that the number of iterations remains
roughly constant and the running time of the solution algorithm scales
roughly linearly with respect to the number of unknowns.

\begin{table}[h!]
  \begin{center}
    \begin{tabular}{|cccc|cc|cc|}
      \hline
      \multicolumn{4}{|c|}{} & \multicolumn{2}{|c|}{Test 1} & \multicolumn{2}{|c|}{Test 2}\\
      \hline
      $\omega/(2\pi)$ & $q$ & $N=n^3$ & $T_{\text{setup}}$ & $N_{\text{iter}}$ & $T_{\text{solve}}$ & $N_{\text{iter}}$ & $T_{\text{solve}}$\\
      \hline
      5 & 8  & $39^3$  & 4.87e+00 & 11 & 4.91e+00 & 11 & 4.96e+00\\
      5 & 16 & $79^3$  & 6.59e+01 & 11 & 4.70e+01 & 12 & 5.55e+01\\
      5 & 32 & $159^3$ & 8.07e+02 & 13 & 5.91e+02 & 13 & 6.31e+02\\
      \hline
    \end{tabular}
  \end{center}
  \caption{Results of velocity field 1 with varying $q$.}
  \label{tbl:3DSPL1}
\end{table}

\begin{table}[h!]
  \begin{center}
    \begin{tabular}{|cccc|cc|cc|}
      \hline
      \multicolumn{4}{|c|}{} & \multicolumn{2}{|c|}{Test 1} & \multicolumn{2}{|c|}{Test 2}\\
      \hline
      $\omega/(2\pi)$ & $q$ & $N=n^3$ & $T_{\text{setup}}$ & $N_{\text{iter}}$ & $T_{\text{solve}}$ & $N_{\text{iter}}$ & $T_{\text{solve}}$\\
      \hline
      5 & 8  & $39^3$  & 4.80e+00 & 12 & 5.36e+00 & 12 & 4.95e+00\\
      5 & 16 & $79^3$  & 6.74e+01 & 13 & 5.53e+01 & 12 & 5.51e+01\\
      5 & 32 & $159^3$ & 8.18e+02 & 14 & 6.48e+02 & 14 & 6.45e+02\\
      \hline
    \end{tabular}
  \end{center}
  \caption{Results of velocity field 2 with varying $q$.}
  \label{tbl:3DSPL2}
\end{table}

\begin{table}[h!]
  \begin{center}
    \begin{tabular}{|cccc|cc|cc|}
      \hline
      \multicolumn{4}{|c|}{} & \multicolumn{2}{|c|}{Test 1} & \multicolumn{2}{|c|}{Test 2}\\
      \hline
      $\omega/(2\pi)$ & $q$ & $N=n^3$ & $T_{\text{setup}}$ & $N_{\text{iter}}$ & $T_{\text{solve}}$ & $N_{\text{iter}}$ & $T_{\text{solve}}$\\
      \hline
      5 & 8  & $39^3$  & 4.82e+00 & 12 & 4.92e+00 & 12 & 5.08e+00\\
      5 & 16 & $79^3$  & 6.77e+01 & 12 & 5.17e+01 & 13 & 6.04e+01\\
      5 & 32 & $159^3$ & 8.16e+02 & 13 & 6.26e+02 & 15 & 7.14e+02\\
      \hline
    \end{tabular}
  \end{center}
  \caption{Results of velocity field 3 with varying $p$.}
  \label{tbl:3DSPL3}
\end{table}

Let us compare these numerical results with the ones from the 3D
results from the previous paper \cite{EngquistYing:10a}. The setup
time $T_{\text{setup}}$ of the current algorithms is much lower: for
the problem of $20$ wavelength across, the current setup time is in
the hundreds of seconds while the setup time in
\cite{EngquistYing:10a} is in the tens of thousands of seconds. This
is mainly due to the fact that our implementation of the multifrontal
algorithm in this paper is more efficient compared to our
implementation of the 2D hierarchical matrix algebra in
\cite{EngquistYing:10a}. The number of iterations $N_{\text{iter}}$ is
about 5 times larger, again due to the fact that the current
algorithms use physical arguments about the Helmholtz equation rather
than direct numerical approximation for $T_m$. Notice that the
solution time $T_{\text{solve}}$ is only about 3 to 4 times larger and
this is due to the efficiency of applying $\wt{T}_m$ using the
multifrontal factorization.

\section{Conclusion and Future Work}
\label{sec:conc}

In this paper, we proposed a new sweeping preconditioner for the
Helmholtz equation in two and three dimensions. Similar to the
previous paper \cite{EngquistYing:10a}, the preconditioner is based on
an approximate block $LDL^t$ factorization that eliminates the
unknowns layer by layer starting from an absorbing layer or boundary
condition. What is new is that the Schur complement matrices of the
block $LDL^t$ factorization are approximated by introducing moving
PMLs in the interior of the domain. In the 2D case, applying these
Schur complement matrices corresponds to solving quasi-1D problems by
an LU factorization with optimal ordering. In the 3D case, applying
these Schur complement matrices corresponds to solving quasi-2D
problems with multifrontal methods. The resulting preconditioner has a
linear application cost and the number of iterations is essentially
independent of the number of unknowns or the frequency when combined
with the GMRES solver.

Some questions remain open. First, we tested the algorithms with the
PML boundary condition as the numerical implementation of the
Sommerfeld condition. Many other boundary conditions are available and
we believe that the current algorithms should work for these boundary
conditions. We presented the algorithms using the simplest central
difference scheme (5 point stencil in 2D and 7 point stencil in
3D). The dispersion relationships of these schemes are rather poor
approximations to the true one. One would like to investigate other
more accurate stencils and other types of discretizations such as
finite element, spectral element, and discontinuous Galerkin.

Parallel processing is necessary for large scale 3D problems. Although
the overall structure of the sweeping preconditioner is sequential by
itself, the calculation of the multifrontal method within each
iteration can be well parallelized. Several efficient implementations
are already available
\cite{AmestoyDuffLExcellentKoster:01,LinYangLuYingE:10} for this
purpose. There is also an alternative to parallelize via a coarse
scale domain decomposition and apply our technique within each
subdomain.

The approach of the current paper is readily applicable to non-uniform
and even adaptive grids. In fact, the same non-uniform or adaptive
grid can be used for the subproblems associated with the moving
PMLs, as long as the grid can resolve the moving PML with
sufficient accuracy. Since the multifrontal methods for non-uniform
and adaptive grids are readily available
\cite{AmestoyDuffLExcellentKoster:01,LinYangMezaLuYingE:10}, it makes
the current approach more flexible compared with the one of the
previous paper \cite{EngquistYing:10a} based on the hierarchical matrix
representation.


\bibliographystyle{abbrv}
\bibliography{ref}

\end{document}